\pdfoutput=1
\RequirePackage{ifpdf}
\ifpdf 
\documentclass[pdftex]{sigma}
\else
\documentclass{sigma}
\fi

\usepackage[all]{xy}
\usepackage{tikz}\usetikzlibrary{matrix,arrows,calc}
\numberwithin{equation}{section}
\newtheorem{Theorem}{Theorem}[section]

\newtheorem{Lemma}[Theorem]{Lemma}
\newtheorem{Proposition}[Theorem]{Proposition}
{\theoremstyle{definition}
\newtheorem{Definition}[Theorem]{Definition}

\newtheorem{Example}[Theorem]{Example}
\newtheorem{Remark}[Theorem]{Remark}
}
\tikzset{node distance=2cm, auto}

\newcommand{\arrows}{\,\lower1pt\hbox{$\longrightarrow$}\hskip-.24in\raise2pt \hbox{$\longrightarrow$}\,}

\newcommand{\xfy}{X \stackrel{f}{\longleftarrow}Y}

\newcommand{\wwxfy}{X \stackrel{[f]}{\longleftarrow}Y}

\newcommand{\ygz}{Y \stackrel{g}{\longleftarrow}Z}

\newcommand{\xfygz}{{X \stackrel{f}{\longleftarrow}Y \stackrel{g}{\longleftarrow}Z}}

\newcommand{\xfzgy}{{X \stackrel{f}{\longleftarrow}Z \stackrel{g}{\longrightarrow}Y}}

\begin{document}

\newcommand{\arXivNumber}{1401.7302}

\allowdisplaybreaks

\renewcommand{\thefootnote}{$\star$}

\renewcommand{\PaperNumber}{100}

\FirstPageHeading

\ShortArticleName{Selective Categories and Linear Canonical Relations}

\ArticleName{Selective Categories and Linear Canonical Relations\footnote{This paper is a~contribution to the Special
Issue on Poisson Geometry in Mathematics and Physics.
The full collection is available at \href{http://www.emis.de/journals/SIGMA/Poisson2014.html}
{http://www.emis.de/journals/SIGMA/Poisson2014.html}}}

\Author{David LI-BLAND and Alan WEINSTEIN}

\AuthorNameForHeading{D.~Li-Bland and A.~Weinstein}

\Address{Department of Mathematics, University of California, Berkeley, CA 94720 USA}
\Email{\href{mailto:libland@math.berkeley.edu}{libland@math.berkeley.edu}, \href{mailto:alanw@math.berkeley.edu}{alanw@math.berkeley.edu}}
\URLaddress{\url{http://math.berkeley.edu/~libland/}, \url{http://math.berkeley.edu/~alanw/}}

\ArticleDates{Received July 22, 2014, in f\/inal form October 20, 2014; Published online October 26, 2014}

\Abstract{A construction of Wehrheim and Woodward circumvents the problem that compositions of smooth canonical
relations are not always smooth, building a~category suitable for functorial quantization.
To apply their construction to more examples, we introduce a~notion of \textit{highly selective category}, in which only
certain morphisms and certain pairs of these morphisms are ``good''.
We then apply this notion to the category $\mathbf{SLREL}$ of linear canonical relations and the result
${\rm WW}(\mathbf{SLREL})$ of our version of the WW construction, identifying the morphisms in the latter with pairs $(L,k)$
consisting of a~linear canonical relation and a~nonnegative integer.
We put a~topology on this category of \textit{indexed linear canonical relations} for which composition is continuous,
unlike the composition in $\mathbf{SLREL}$ itself.
Subsequent papers will consider this category from the viewpoint of derived geometry and will concern quantum
counterparts.}


\Classification{53D50; 18F99; 81S10}

\renewcommand{\thefootnote}{\arabic{footnote}}
\setcounter{footnote}{0}

\section{Introduction}

In mathematical formulations of the quantization problem, one seeks correspondences or other connections between
mathematical settings for classical and quantum mechanics.
For instance, on the classical (or particle) side one may consider symplectic manifolds as phase spaces and
symplectomorphisms as symmetries and time evolution operators, while on the quantum (or wave) side one has Hilbert
spaces and unitary operators.
Though no universal correspondence of this type can be perfectly compatible with the group structures on the operators,
compatible correspondences can be def\/ined on many interesting subgroups of symplectomorphism groups, at worst up to
extensions of symplectomorphism groups associated with metaplectic structures.
In particular, the machinery of geometric quantization provides such a~correspondence for a~double covering of the group
${\rm Sp}(V)$ of linear automorphisms of a~f\/inite-dimensional symplectic vector space~$V$.\footnote{In this paper, we will not
be discussing the important subject of deformation quantization, in which the connection between classical and quantum
mechanics is realized by deformations of algebras of observables.}

Many questions related to quantization go beyond automorphisms, though, to more general morphisms between symplectic
manifolds and between Hilbert spaces, such as symplectic reductions on the classical side and projections on the quantum
side.
Microlocal analysis, beginning with work of Maslov~\cite{ma:theory} and H\"ormander~\cite{ho:fourier} suggests that the
appropriate classical morphisms are canonical relations, i.e.~Lagrangian submanifolds of products.

It has been advocated for some time (for geometric quantization in~\cite{we:symplecticgeometry}
and~\cite{we:symplecticcategory}, and for deformation quantization in~\cite{la:quantization}) that good solutions of
quantization problems should involve functors between categories on the classical and quantum sides.
For geometric quantization of symplectic vector spaces, Guillemin and Sternberg~\cite{gu-st:problems} constructed
a~correspondence where the classical morphisms are linear canonical relations and the quantum morphisms are certain
unbounded operators with distributional kernels.
(See also~\cite{tu-za:category}, where these kernels are called ``Fresnel kernels''.) Their correspondence fails,
though, when the composition of a~pair of canonical relations fails to satisfy a~transversality condition, in which case
the composition of the operators no longer has a~distributional kernel, and may even have its domain of def\/inition
reduced to zero.
Beyond the linear setting, worse things yet occur in the absence of transversality (and embedding) assumptions: the
composition of smooth canonical relations can take us outside what might have been thought to form the appropriate
category.

One aim of this paper is to redef\/ine the category of linear symplectic relations in such a~way that composition becomes
continuous, and the category is more amenable to quantization.
First though, we set up a~general approach to dealing with ``improper'' compositions in categories, following
a~construction introduced by Wehrheim and Woodward~\cite{we-wo:functoriality} (to be referred to as WW) for dealing with
singular compositions of smooth canonical relations in symplectic topology.

Our basic idea is to select, within an underlying category ${\mathbf C}$, a~collection of nice morphisms, which we call
\textit{suave}, and a~collection of nice composable pairs of suave morphisms, which we call \textit{congenial}, such
that every pair including an identity morphism is congenial, and such that the composition of a~congenial pair is always
suave.
From such a~\textit{selective category} ${\mathbf C}$, we form the category ${\rm WW}({\mathbf C})$ generated by the suave
morphisms, with relations given by the congenial compositions.
For certain purposes, we also distinguish subcategories of the suave morphisms whose members we call \textit{reductions}
and \textit{coreductions}.
After the imposition of some further axioms, notably the requirement that each suave morphism be the composition of
a~reduction with a~coreduction, we have the def\/inition of a~\textit{highly selective category}.

We extend to all highly selective categories ${\mathbf C}$ the main result of~\cite{we:note} to the ef\/fect that, in the
special case where ${\mathbf C}$ is a~certain highly selective category of relations between symplectic manifolds, with
suave morphisms the smooth canonical relations, every morphism\footnote{We will usually denote morphisms in categories
by arrows $X\leftarrow Y$ from right to left; the notation ${\rm Hom}(X,Y)$ will therefore denote morphisms {\em to}~$X$
{\em from}~$Y$.} $X \leftarrow Y$ in ${\rm WW}({\mathbf C})$ may be represented by a~composition (not necessarily congenial)
$X \twoheadleftarrow Q \leftarrowtail Y$ of just two suave morphisms in ${\cal C}$, where the decorations on the arrows
mean that $X \twoheadleftarrow Q$ is a~reduction and $Q \leftarrowtail Y$ is a~coreduction.

We pay special attention to rigid monoidal structures, with monoidal products denoted $X\otimes Y$, dual objects
${\overline{X}}$, and unit object $\mathbf{1}$.
When a~rigid monoidal structure is compatible with a~selective structure on ${\mathbf C}$, it extends to ${\rm WW}({\mathbf
C})$.
In a~short digression, we present the formalism of string diagrams for portraying computations in rigid monoidal
categories.

In any rigid monoidal category, the morphisms from $\mathbf{1}$ play a~special role.
In fact, ${\rm Hom}(X,Y)$ can be identif\/ied with ${\rm Hom}(X \otimes {\overline{Y}}, \mathbf{1})$, with each morphism
$X\leftarrow Y$ corresponding to its ``graph'' $X\times {\overline{Y}}\leftarrow\mathbf{1}$.
The endomorphisms of $\mathbf{1}$ function as a~``coef\/f\/icient ring''.
These morphisms are represented in the symplectic situation by intersecting pairs of Lagrangian submanifolds, and on the
quantum side by linear pairings of vectors, in each case modulo an equivalence relation whose study we merely begin
here.
In ${\rm WW}({\mathbf C})$, therefore, each morphism may be represented as a~composition $X\otimes {\overline{Y}}
\twoheadleftarrow Q \leftarrowtail \mathbf{1}$.

In categories of linear operators, the morphisms $X\otimes {\overline{Y}}\leftarrow\mathbf{1}$ are sometimes known as
kernels of the corresponding morphisms $X\leftarrow Y$, so the WW-morphisms may be referred to as \textit{hyperkernels}.
Similarly, the morphisms in WW-categories built from categories of relations may be called \textit{hyperrelations}.
We will sometimes refer to a~diagram $X\otimes {\overline{Y}} \twoheadleftarrow Q \leftarrowtail \mathbf{1}$ in any
category as a~\textit{hypermorphism} to~$X$ from~$Y$.
The composition of hypermorphisms is essentially a~monoidal product, which is always def\/ined without any extra
assumptions, though one must remember the equivalence relation, for instance to see why an identity morphism really is
an identity.

In the symplectic case, when~$Y$ is a~point, the Lagrangian submanifolds of~$X$ should be replaced by ``Lagrangian
hypersubmanifolds'', which are Lagrangian submanifolds in symplectic manifolds of which~$X$ is a~symplectic
reduction\footnote{This notion of Lagrangian hypersubmanifold was originally inspired by the special case of generating
functions for Lagrangian submanifolds of cotangent bundles.
Here,~$Z$ is a~cotangent bundle $T^*M$,~$Q$ is a~cotangent bundle~$T^*B$, the reduction morphism $Z\leftarrow Q$ is the
one associated to (i.e.~the cotangent lift of) a~surmersion of manifolds~$M\leftarrow B$, and~$L$ is the image of the
dif\/ferential of a~function~$S\in {C^\infty}(B)$.
$L$ is transversal to the domain of the reduction relation (the
conormal bundle to the f\/ibres of $M\leftarrow B$) just when~$S$ is a~Morse family over~$M$; in this case, the reduction
of~$L$ in $T^*M$ is the (immersed, in general) Lagrangian submanifold for which~$S$ is a~generating function.}.

After developing the general WW-theory, we consider situations where the morphisms of the orginal category are
set-theoretic relations.
There are many interesting but dif\/f\/icult problems concerning the classif\/ication of WW-morphisms built from relations,
related to a~notion of ``trajectory''.
In particular, we look at relations between symplectic manifolds, with the suave ones being those given by closed
Lagrangian submanifolds of products.
(This example was already treated in~\cite{za:quantum}.)

\looseness=1
After some brief remarks on classif\/ication of hypermorphisms $\mathbf{1} \leftarrow \mathbf{1}$ in the symplectic
mani\-fold setting, we concentrate for the remainder of the paper on {\em linear} canonical (hyper)re\-la\-tions.

Although the linear canonical relations between f\/inite-dimensional vector spaces already form a~category
$\mathbf{SLREL}$ (see~\cite{be-tu:relazionisimplettiche} for an early discussion of this category) before the WW
construction is applied, this category has the defect that the composition operation is not continuous in the usual
topology on the morphism spaces, which may be identif\/ied with Lagrangian Grassmannians.
For this reason, and for purposes of quantization, we introduce in $\mathbf{SLREL}$ the selective structure in which all
morphisms are suave, but only transversal compositions are congenial.
Here, we can completely classify the WW-morphisms, identifying them with pairs $(L,k)$, where~$L$ is a~canonical
relation and~$k$ is a~nonnegative integer.
We call these \textit{indexed canonical relations}.
Thus, the WW endomorphisms of the unit object (a zero vector space) are just nonnegative integers, with composition
given by addition.
The composition operation is now continuous in a~topology which we will describe in detail.
We call it the \textit{Sabot topology} since its def\/inition is inspired by the description in~\cite{sa:electrical} of
the closure of the graph of a~(discontinuous) symplectic reduction map between Lagrangian Grassmannians.
In fact, using Sabot's result and ref\/inements thereof, we show that the Sabot topology is essentially obtained from
a~quotient construction, via symplectic reduction, from the usual topology on the Lagrangian Grassmannian.

We turn now to quantization, building on earlier work of Guillemin and Sternberg~\cite{gu-st:problems}, Tul\-czyjew and
Zakrezewski~\cite{tu-za:category}, and Benenti~\cite{be:category}.
These authors quantize linear canonical relations by (possibly unbounded) operators\footnote{Actually, each operator is
def\/ined only up to a~multiplicative constant.
Choosing the constant requires some extra ``amplitude'' and ``metaplectic'' data on the symplectic side.} whose Schwartz
kernels are, after a~suitable splitting of variables, built from Dirac delta ``functions'' combined with exponentials of
imaginary quadratic polynomials.
They show that the quantization of the composition of two relations is the product of the corresponding operators, as
long as the classical composition is congenial in the sense that a~transversality condition is satisf\/ied.
When the condition fails, the product of operators does not have a~Schwartz kernel because its domain has become too
small, so the quantum composition does not even exist.

The underlying problem is that the operators with Schwartz kernels fail to constitute the morphisms of a~category.
One of the simplest cases of this failure comes down to the well-known dif\/f\/iculty of multiplying the delta function
$\delta(x)$ in a~single variable by itself, or, more or less equivalently, evaluating $\delta(x)$ at its singular
point~$0$.
This issue was treated symplectically to some (insuf\/f\/icient) extent in~\cite{we:maslovcycle} as well as, earlier,~by
Sabot~\cite{sa:electrical} from a~completely dif\/ferent point of view.

In the companion paper~\cite{jo-li-we:linearsymplectic}, we will describe a~``quantum WW construction,'' based on an
idea described in a~letter to the second author from Shlomo Sternberg (written in 1982 on the shore of the Red Sea but
only recently excavated from the f\/ile drawer in which it was buried), who attributed it to a~seminar talk by Ofer
Gabber.
The basic idea, quite natural to experts in D-module theory, is to quantize Lagrangian subspaces not by functions, or
even by distributions, but by modules over the Weyl algebras associated to symplectic vector spaces.
For canonical relations, these modules are bimodules, and the composition is given by a~tensor product.
Even here, functoriality fails for nontransversal compositions, and it is necessary to take a~further step, from modules
to complexes of modules, and from tensor products to torsion products.
This cohomological approach, going back at least to Serre's fundamental work~\cite{se:algebre} on intersection theory,
is a~simple instance of what is now sometimes known as ``derived noncommutative geometry''.
One reference for this approach is~\cite{ka-sc:sheaves}.

In~\cite{jo-li-we:linearsymplectic}, before developing in detail the idea of Gabber described above, we will introduce
a~parallel structure on the classical side.
Namely, we will consider symplectic vector spaces and their Lagrangian subspaces as a~special case (concentrated in
degree 0) of ``symplectic complexes'' of vector spaces and ``Lagrangian morphisms'' from other complexes to these.
The composition of morphisms is def\/ined via derived tensor products so that, when pairs in degree zero are composed, the
composition may live in other degrees, depending upon the extent of nontransversality of the composition.
We can then connect this category with the WW-category and the more concrete category of indexed canonical relations.
In particular, the Sabot topology will be seen as the natural topology on the coarse moduli space of a~stack of
Lagrangian subspaces in derived symplectic geometry.

It may still be possible to f\/ind more concrete approaches to the composition problem on the quantum side, using some
combination of the theories of partial inner product spaces~\cite{an-tr:partial} and Columbeau's generalized
functions~\cite{co:new}.
We leave this for future work.

\section{Selective categories}

The embedding by Wehrheim and Woodward~\cite{we-wo:functoriality} of the canonical relations between symplectic
manifolds into an actual category depends only on some extra structure on the category of {\em all} relations between
symplectic manifolds.
We will formalize this kind of structure below and give many examples, some of which may already be found
in~\cite{we:note}.

\begin{Definition}
\label{dfn-selective}
A~\textit{selective category} is a~category together with a~distinguished class of morphisms, called \textit{suave}, and
a~class of composable pairs of suave morphisms called \textit{congenial pairs}, such that:
\begin{enumerate}\itemsep=0pt
\item Any identity morphism is suave.
\item If~$f$ and~$g$ are suave, $(f,g)$ is composable, and~$f$ or~$g$ is an identity morphism, then $(f,g)$ is
congenial.
\item If $(f,g)$ is congenial, then $fg$ is suave.
\item If~$f$ is a~suave isomorphism, its inverse $f^{-1}$ is suave as well, and the pairs $(f,f^{-1})$ and $(f^{-1},f)$
are both congenial.
\item If $(f,g,h)$ is a~composable triple, then $(f,g)$ and $(fg,h)$ are congenial if and only if $(g,h)$ and $(f,gh)$
are.
When these conditions hold, we call $(f,g,h)$ a~\textit{congenial triple}.
\end{enumerate}
\end{Definition}

When $(f,g)$ is congenial, we will sometimes refer to the expression $fg$ as a~\textit{congenial composition}.

\begin{Remark}
Just as the associative law implies the generalized associative law for removal of parentheses in products of four or
more elements, so the last condition in the def\/inition above implies that the same proof idea can be used to show that,
if a~product of suave isomorphisms can be computed congenially with respect to one parenthesization, the same is true
for any parenthesization.
Thus, there is a~well-def\/ined notion of congenial~$n$-tuple for any~$n$.
\end{Remark}
Another consequence of the ``associativity of congeniality'' is the following.
\begin{Proposition}
\label{prop-compwithisomo}
A~pair $(f,g)$ of suave morphisms is congenial whenever either factor is an isomorphism.
\end{Proposition}
\begin{proof}
Suppose that~$f$ is a~suave isomorphism.
(The case when~$g$ is an isomorphism is similar.) Then $f^{-1}$ is suave as well, and the pairs $(f^{-1},f)$ and
$(f^{-1}f,g)$ are congenial, so $(f^{-1},f,g)$ is a~congenial triple; hence $(f,g)$ is congenial.
\end{proof}

\begin{Example}
\label{ex-selectivesuave}
We denote by $\mathbf{REL}$ the category whose objects are sets and whose morphisms are relations; i.e., the morphisms
in $\mathbf{REL}(X,Y)$ are the subsets of $X\times Y$.

We will use the selective structure in $\mathbf{REL}$ in which all morphisms are suave, but only \textit{monic} pairs
are congenial.
These are diagrams $\xfygz$ for which, whenever $(x,y)$ and $(x,y')$ belong to~$f$, and $(y,z)$ and $(y',z)$ belong
to~$g$, then $y=y'$.
In other words, $(x,z)$ can belong to $f\circ g$ ``in at most one way''.
The only condition in Def\/inition~\ref{dfn-selective} which takes a~bit of work to check is the last one.
For this, it suf\/f\/ices to observe that congenial triples can be def\/ined directly: $(f,g,h)$ is congenial if and only if,
whenever $(x,w)$ belongs to $fgh$, there is exactly one pair $(y,z)$ for which $(x,y) \in f$, $(y,z) \in g$, and $(z,w)
\in h$.

We will be concerned later with categories of relations between sets with additional structure.
For instance, $\mathbf{MREL}$ and $\mathbf{SREL}$ will denote the categories whose objects are smooth manifolds and
symplectic manifolds respectively, but with all set-theoretic relations as morphisms.
Although the forgetful functors
\begin{gather*}
\mathbf{REL} \leftarrow \mathbf{MREL} \leftarrow \mathbf{SREL}
\end{gather*}
are fully faithful, we are not just dealing with subcategories here, since any uncountable set carries many dif\/ferent
smooth and symplectic structures.

In $\mathbf{MREL}$, we will use the selective structure in which the suave morphisms are the smooth relations
(i.e.~closed submanifolds of products).
In $\mathbf{SREL}$, the suave morphisms will be the canonical relations, i.e.~those $f\in \mathbf{MREL}(X,Y)$ which are
closed Lagrangian submanifolds of $X\times {\overline{Y}}$.\footnote{If~$Y$ is a~symplectic manifold, ${\overline{Y}}$
is the same manifold with the symplectic structure multiplied by $-1$.
If~$Y$ is just a~manifold, ${\overline{Y}}$ is equal to~$Y$.} In both cases, the congenial pairs will be the strongly
transversal pairs as def\/ined in~\cite[Def\/inition~3.3]{we:note}; i.e.~they are the $(f,g)$ for which $(f\times g) \cap
(X \times \Delta_Y\times {\overline{Z}})$ is a~transversal intersection, and the projection of this intersection to $X
\times {\overline{Z}}$ is a~proper embedding.
(The second condition is a~strong version of monicity.) In the symplectic case, transversality is equivalent to the
local version of monicity in which the projection of $(f\times g) \cap (X \times \Delta_Y\times Z)$ to $X\times
{\overline{Z}}$ is an immersion.
\end{Example}

\begin{Definition}
A~\textit{selective functor} between selective categories is one which takes congenial pairs to congenial pairs.
\end{Definition}

The forgetful functors above are both selective.

Composition with identity morphisms shows that a~selective functor takes suave morphisms to suave morphisms.

In many selective categories of interest, the collection of suave morphisms is not closed under composition, so these
are not the morphisms of a~category.
The Wehrheim--Woodward construction circumvents this problem (and others) by faithfully embedding the suave morphisms in
${\mathbf C}$ as morphisms in a~category ${\rm WW}({\mathbf C})$ in which the composition of any {\em congenial} pair remains
the same as the composition in ${\mathbf C}$.

The construction of ${\rm WW}({\mathbf C})$ begins with a~category of ``paths'' in ${\mathbf C}$ (not using the selective
structure).
If we think of a~category as a~directed graph with the objects as vertices and morphisms as edges, these are paths in
the usual sense, up to ``weakly monotonic'' reparametrization, as is made precise below.

\begin{Definition}
The \textit{support} of an inf\/inite composable sequence
\begin{gather*}
f=(\ldots, f_{-1},f_0,f_1,\ldots)
\end{gather*}
of morphisms in a~category ${\mathbf C}$ is the set of integers~$j$ for which $f_j$ is not an identity morphism.
A~\textit{path in} ${\mathbf C}$ is an inf\/inite composable sequence of f\/inite support.
The target and source of $f_j$ for all suf\/f\/iciently large negative~$j$ is thus a~f\/ixed object~$X$ and, for all
suf\/f\/iciently large positive~$j$, a~f\/ixed object~$Y$.
We call~$X$ the \textit{target} and~$Y$ the \textit{source} of the path~$f$.

Two paths will be considered as equivalent if one may be obtained from the other by inserting and removing f\/initely many
identity morphisms.
This does not change the target or source.
The set of equivalence classes is the \textit{path category} $P({\mathbf C})$.
We will denote the equivalence class of $f=(\ldots, f_{-1},f_0,f_1,\ldots)$ by $ \langle \ldots, f_{-1},f_0,f_1,\ldots
\rangle$, or simply $\langle f \rangle$.
We will also use the notation $\langle f_r,\ldots,f_s \rangle $ when the support of the sequence~$f$ is contained in the
interval $[r,s]$.

To compose $ \langle f\rangle \in P({\mathbf C})(X,Y)$ and $\langle g \rangle \in P({\mathbf C})(Y,Z)$, choose
representative sequences, remove all but f\/initely many consecutive copies of $1_Y$ from the positive end of the f\/irst
sequence and from the negative end of the second, and then concatenate the truncated sequences.

The identity morphism in $P({\mathbf C})$ of any object~$X$ is (represented by) the constant sequence with all entries
equal to $1_X$.

If ${\mathbf C}$ is a~selective category, we denote the subcategory of $P({\mathbf C})$ consisting of equivalence
classes of paths of suave morphisms by $P_s({\mathbf C})$.
\end{Definition}

\begin{Remark}
\label{rmk-minimal}
We could have def\/ined paths as (equivalence classes of) f\/inite sequences, but the inf\/inite version is more convenient
when it comes to def\/ining rigid monoidal structures.
Still, it is convenient to have f\/inite representations of paths.

Every morphism $\langle f \rangle$ in $P({\mathbf C})$ has a~unique ``minimal'' representative for which $f_i$ is an
identity morphism for all $i\leq 0$ and for which there are no identity morphisms between two nonidentity morphisms.
It may therefore be denoted in the form $\langle f_1,\ldots,f_n\rangle$, where none of the $f_i$ is an identity morphism
(except when $\langle f \rangle$ is itself an identity morphism).
We will frequently use this notation, often without the elimination of identity morphisms.
We will also use the f\/inite representation $(f_1,\ldots,f_n)$ for the paths themselves.
\end{Remark}

\begin{Remark}
\label{rmk-shifting}
A~useful way to carry out the composition of two sequences is to shift the f\/irst one (which does not change its
equivalence class) so that its support is contained in the negative integers, and to shift the second so that its
support is contained in the positive integers.
The composition is then represented by the sequence whose value at~$j$ is $f_j$ for $j\leq 0$ and $g_j$ for $j\geq 0$.

One may use a~similar idea to verify associativity of composition; given three sequences, shift them so that their
supports are contained in disjoint, successive intervals of integers.
\end{Remark}

We leave to the reader the proof of the following result.
\begin{Proposition}
\label{prop-pathcomposition}
There is a~unique functor ${\mathbf C} \stackrel{c'}{\leftarrow} P({\mathbf C})$ which is the identity on objects and
which takes each morphism $ \langle f_{1},\ldots,f_n\rangle $ to the composition $ f_{1} \cdots f_n $
in ${\mathbf C}$.
\end{Proposition}

We now def\/ine the Wehrheim--Woodward
category ${\rm WW}({\mathbf C})$ by permitting the composition of congenial pairs
contained in paths.

\begin{Definition}
Let ${\mathbf C}$ be a~selective category.
The Wehrheim--Woodward category ${\rm WW}({\mathbf C})$ is the quotient category obtained from $P_s({\mathbf C})$ by the
smallest equivalence relation for which two paths are equivalent if a~sequence representing one is obtained from
a~sequence representing the other by replacing successive entries forming a~congenial pair $(p,q)$ by the composition
$pq$ preceded or followed by an identity morphism.
The equivalence class in ${\rm WW}({\mathbf C})$ of $\langle f \rangle \in P_s({\mathbf C})$ will be denoted by $[f]$.
\end{Definition}

\begin{Remark}
As was the case for paths (see Remark~\ref{rmk-minimal}), every morphism in ${\rm WW}({\mathbf C})$ has a~representation of
the form $[f_1,\ldots,f_n]$.
\end{Remark}

\begin{Remark}
The identity morphism preceding or following $pq$ is not essential, since it can be removed, but it is sometimes helpful
to keep the remaining entries in the sequence unshifted.
\end{Remark}

Here is another useful fact.

\begin{Proposition}
Let $(f'_1,\ldots,f'_n)$ and $(f''_1,\ldots,f''_n)$ be paths representing morphisms to~$X$ from~$Y$ in ${\rm WW}({\mathbf
C})$, with $X'_{j-1}$ and $X'_j$ $($resp.\ $X''_{j-1}$ and $X''_j)$ the target and source objects of $f'_j$ $($resp.\ $f''_j)$.
$($In particular, $X'_0 = X= X''_0$ and $X'_n = Y = X''_n.)$ If the two paths are isomorphic in the sense that there exist
suave isomorphisms $X'_j \stackrel{\phi_j}{\leftarrow} X''_j$ such that $\phi_0$ and $\phi_n$ are identity morphisms,
and $\phi_{j-1} f''_j = f'_j \phi_j$ for all~$j$, then the $\rm WW$-morphisms $[f'_1,\ldots,f'_n]$ and $[f''_1,\ldots,f''_n]$
are equal.
\end{Proposition}

\begin{proof}
Since composition with a~suave isomorphism is congenial, we have $[f''_j] = [\phi_{j-1}^{-1} f'_j \phi_j]$ for
$j=1,\ldots,n$.
Multiplying these~$n$ equalities gives the desired result.
\end{proof}

There is a~map~$\iota$ from the suave morphisms in ${\mathbf C}$ to ${\rm WW}({\mathbf C})$ which takes each suave
morphism~$f$ to the equivalence class of sequences containing one entry equal to~$f$ and all the others equal to
identity morphisms.
${\rm WW}({\mathbf C})$ is then characterized by the following universal property, whose proof we omit. 

\begin{Proposition}
The composition functor $c'$ above descends to a~functor ${\mathbf C} \stackrel{c}{\leftarrow} {\rm WW}({\mathbf C})$, namely
$c([f_{1},\ldots,f_n]) = f_{1}\cdots f_n$. More generally, ${\rm WW}({\mathbf C})$ has the universal property that any map from
the suave morphisms of ${\mathbf C}$ to a~category ${\mathbf B}$ which takes units to units and which takes congenial
compositions to compositions in ${\mathbf B}$ is of the form $b\circ \iota$ for a~unique functor ${\rm WW}({\mathbf
C})\stackrel{b}{\rightarrow}{\mathbf B}$.
\end{Proposition}

Applying this proposition to the inclusion of the congenial pairs in ${\mathbf C}$, we may conclude that~$\iota$ is
injective, i.e.~distinct suave morphisms cannot become equal when considered as WW-morphisms via~$\iota$.

We will refer to $c([f_{1},\ldots,f_n])$ as the \textit{shadow} of $[f_{1},\ldots,f_n]$.

\begin{Remark}
Any selective functor between selective categories induces a~functor between their Wehrheim--Woodward categories.
In particular, if a~selective category ${\mathbf C}$ admits a~transpose operation, i.e.~an involutive contravariant
endofunctor $f\mapsto f^t$ which f\/ixes objects and takes congenial pairs to congenial pairs, then this operation extends
to ${\rm WW}({\mathbf C})$, where the transpose of a~sequence is the same sequence with its order reversed and each entry
replaced by its transpose.
\end{Remark}

\begin{Remark}
If the suave morphisms in ${\mathbf C}$ form a~subcategory ${\mathbf C} '$ and we declare every composable pair to be
congenial, then ${\mathbf C} ' \stackrel{c}{\leftarrow} {\rm WW}({\mathbf C})$ is an isomorphism of categories.
\end{Remark}

\section{Highly selective categories}

\begin{Definition}
\label{dfn-highly}
A~\textit{highly selective category} is a~selective category provided with two subcategories of suave morphisms called
\textit{reductions} and \textit{coreductions} such that:
\begin{enumerate}\itemsep=0pt
\item Any suave isomorphism is a~coreduction and a~reduction.
\item If $(f,g)$ is a~composable pair of suave morphisms, and if~$f$ is a~coreduction or~$g$ is a~reduction, then
$(f,g)$ is congenial.
\item Any suave morphism~$f$ may be factored as $gh$, where~$g$ is a~reduction,~$h$ is a~coreduction, and $(g,h)$ is
congenial.
\end{enumerate}
\end{Definition}

It follows from the injectivity of~$\iota$ that the subcategories of reductions and coreductions in ${\mathbf C}$ are
mapped isomorphically to subcategories of ${\rm WW}(\mathbf C)$, which we will again refer to as reductions and coreductions.
Since all identity morphisms are suave, these subcategories are wide; i.e., they contain all the objects.

We will indicate that a~morphism is special by decorating its arrow: $X \twoheadleftarrow Q$ is a~reduction, and $Q
\leftarrowtail Y$ is a~coreduction.
\begin{Example}
\label{ex-relhighlyselective}
In the selective category of relations in Example~\ref{ex-selectivesuave}, with the congenial pairs the monic ones, we
def\/ine a~highly selective structure by letting the reductions be the single-valued surjective relations and the
coreductions the everywhere-def\/ined injective ones.
The required factorization of a~morphism $\xfy$ may be constructed as follows.
Take the reduction $X\stackrel{g}{\twoheadleftarrow} X\times Y\times Y$ to be $\{(x,(x',y',y''))| (x=x',y'= y'')\}$ and
the coreduction $ X\times Y\times Y \stackrel{h}{\leftarrowtail} Y $ to be $\{((x',y',y''),y) | (x',y')\in f, y''=
y)\}$.

The same ideas work in $\mathbf{MREL}$ and $\mathbf{SREL}$.
In either case, we can def\/ine a~reduction to be a~reduction in the sense above, together with the condition that the
domain be a~closed submanifold of the source, and the map from this domain to the target a~surjective submersion.
We can even require that this map be a~locally trivial f\/ibration to get a~``f\/iner'' highly selective structure.
Coreductions are just the transposes of reductions.

The construction here extends to any rigid monoidal category, as we will see in Section~\ref{sec-rigidmonoidal}.
\end{Example}

The next result was proven in~\cite{we:note} for the special case of $\mathbf{SREL}$.
The proof in the general case is identical, but we repeat it here for readers' convenience.

\begin{Theorem}
\label{thm-two morphisms}
Let ${\mathbf C}$ be a~highly selective category, and let $[f_1,\ldots,f_n]$ be a~morphism in ${\rm WW}({\mathbf C})$.
Then there exist an object~$Q$ and morphisms $A \in {\rm Hom}(X_0,Q)$ and $B \in {\rm Hom}(Q,X_n)$ in ${\mathbf C}$ such
that~$A$ is a~reduction,~$B$ is a~coreduction, and $[f_1,\ldots,f_n]=[A,B]$ in ${\rm WW}({\mathbf C})$.
\end{Theorem}

\begin{proof}
We illustrate the proof with diagrams for the case $n=4$, which is completely representative of the general case.

First, we write $[f_{1},f_{2},f_{3},f_{4}]$ as a~composition:
\begin{gather*}
\xymatrix{
X_0  && X_1  \ar[ll]_{f_{1}}  &&  X_2
\ar[ll]_{f_{2}} &&    X_3 \ar[ll]_{f_{3}}
&&  X_4 \ar[ll]_{f_{4}}
}
\end{gather*}

Since ${\mathbf C}$ is highly selective, we may factor each arrow as the congenial composition of a~reduction and
a~coreduction.
The top row of the next diagram is then equivalent to the zigzag line below it
\begin{gather*}
\xymatrix{
X_0  && X_1  \ar@{>->}[dl] \ar[ll]_{f_{1}}  &&  X_2
\ar@{>->}[dl]\ar[ll]_{f_{2}} &&  \ar@{>->}[dl]  X_3 \ar[ll]_{f_{3}}
&&  X_4 \ar@{>->}[dl]\ar[ll]_{f_{4}}
\\
  & X_{01}  \ar@{->>}[ul]&& X_{12}  \ar@{->>}[ul]&&
            X_{23} \ar@{->>}[ul]
      && X_{34} \ar@{->>}[ul]\\
 }
\end{gather*}

Next, we compose pairs of diagonal arrows to produce the bottom row below
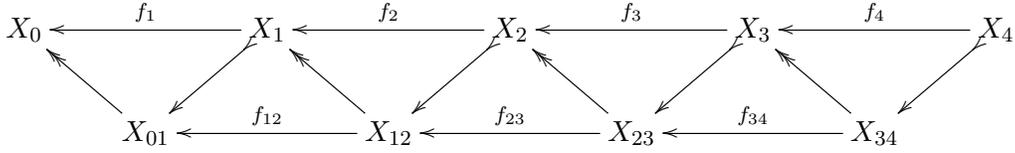
\vspace{-2mm}
\begin{gather*}
\xymatrix{
X_0  && X_1  \ar@{>->}[dl]\ar[ll]_{f_{1}}  &&  X_2
\ar@{>->}[dl]\ar[ll]_{f_{2}} &&  \ar@{>->}[dl]  X_3 \ar[ll]_{f_{3}}
&&  X_4 \ar@{>->}[dl]\ar[ll]_{f_{4}}
\\
  & X_{01}  \ar@{->>}[ul]&& X_{12}  \ar@{->>}[ul]\ar[ll] _{f_{12}}&&
            X_{23} \ar@{->>}[ul]\ar[ll] _{f_{23}}
      && X_{34} \ar@{->>}[ul]\ar[ll]  _{f_{34}} \\
 }
\end{gather*}

\vspace{-2mm}

Each of these compositions is congenial, even ``doubly so'', thanks to the decorations on the arrows identifying them as
reductions and coreductions.
It follows that the original composition on the top row is equivalent to the composition of the bottom row with the
outer diagonal edges.

We repeat the process to obtain another row
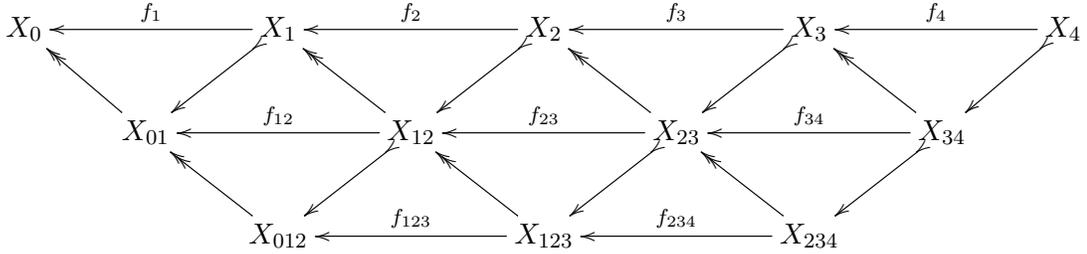
\vspace{-2mm}
\begin{gather*}
\xymatrix{
X_0  && X_1  \ar@{>->}[dl]\ar[ll]_{f_{1}}  &&  X_2
\ar@{>->}[dl]\ar[ll]_{f_{2}} &&  \ar@{>->}[dl]  X_3 \ar[ll]_{f_{3}}
&&  X_4 \ar@{>->}[dl]\ar[ll]_{f_{4}}
\\
  & X_{01}  \ar@{->>}[ul]&& X_{12} \ar@{>->}[dl] \ar@{->>}[ul]\ar[ll] _{f_{12}}&&
            X_{23}  \ar@{>->}[dl]\ar@{->>}[ul]\ar[ll] _{f_{23}}
      && X_{34} \ar@{->>}[ul]\ar@{>->}[dl]\ar[ll]  _{f_{34}}\\
   && X_{012}  \ar@{->>}[ul]&& X_{123} \ar@{->>}[ul]\ar[ll] _{f_{123}} &&
   X_{234} \ar@{->>}[ul]\ar[ll] _{f_{234}}
}
\end{gather*}

\vspace{-2mm}

Repeating two more times, we arrive at a~triangle, in which the top row is equivalent to the composition of the arrows
on the other two sides
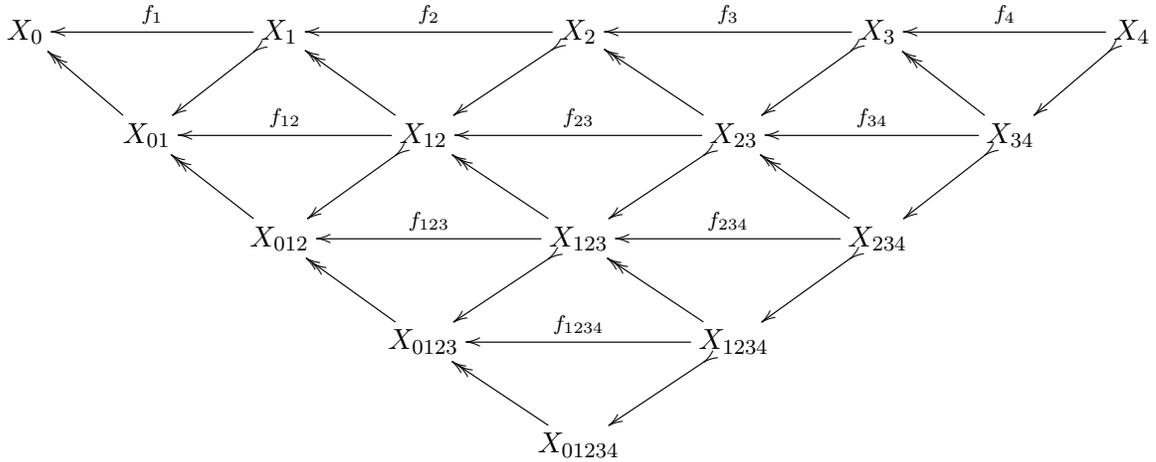
\vspace{-2mm}
\begin{gather*}\hspace*{-10.6pt}
\xymatrix{
X_0  && X_1  \ar@{>->}[dl]\ar[ll]_{f_{1}}  &&  X_2
\ar@{>->}[dl]\ar[ll]_{f_{2}} &&  \ar@{>->}[dl]  X_3 \ar[ll]_{f_{3}}
&&  X_4 \ar@{>->}[dl]\ar[ll]_{f_{4}}
\\
  & X_{01}  \ar@{->>}[ul]&& X_{12} \ar@{>->}[dl] \ar@{->>}[ul]\ar[ll] _{f_{12}}&&
            X_{23}  \ar@{>->}[dl]\ar@{->>}[ul]\ar[ll] _{f_{23}}
      && X_{34} \ar@{->>}[ul]\ar@{>->}[dl]\ar[ll]  _{f_{34}}\\
   && X_{012}  \ar@{->>}[ul]&& X_{123} \ar@{>->}[dl] \ar@{->>}[ul]\ar[ll] _{f_{123}} &&
   X_{234} \ar@{->>}[ul]\ar @{>->}[dl]\ar[ll] _{f_{234}}
\\
     &&&  X_{0123} \ar@{->>}[ul]&& X_{1234}\ar @{>->}[dl] \ar@{->>}[ul]\ar[ll] _{f_{1234}}\\
      &&&&
X_{01234}  \ar @{->>}[ul]
}
\end{gather*}

\vspace{-2mm}

Finally, we observe that all the arrows going up the left-hand side are reductions, so we may compose them all to
produce a~single reduction~$A$.
Similarly, the coreductions going down on the right yield a~coreduction~$B$
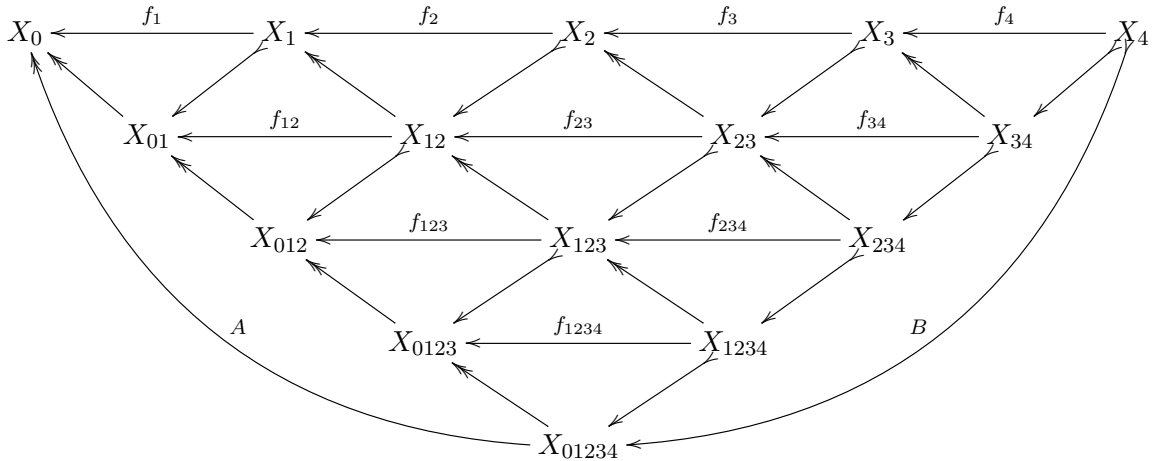
\vspace{-2mm}
\begin{gather*}\hspace*{-10.6pt}
\xymatrix{
X_0  && X_1  \ar@{>->}[dl]\ar[ll]_{f_{1}}  &&  X_2
\ar@{>->}[dl]\ar[ll]_{f_{2}} &&  \ar@{>->}[dl]  X_3 \ar[ll]_{f_{3}}
&&  X_4 \ar@{>->}[dl]\ar[ll]_{f_{4}}
\ar@/^4pc/@{>->}[ddddllll]_B
\\
  & X_{01}  \ar@{->>}[ul]&& X_{12} \ar@{>->}[dl] \ar@{->>}[ul]\ar[ll] _{f_{12}}&&
            X_{23}  \ar@{>->}[dl]\ar@{->>}[ul]\ar[ll] _{f_{23}}
      && X_{34} \ar@{->>}[ul]\ar@{>->}[dl]\ar[ll]  _{f_{34}}\\
   && X_{012}  \ar@{->>}[ul]&& X_{123} \ar@{>->}[dl] \ar@{->>}[ul]\ar[ll] _{f_{123}} &&
   X_{234} \ar@{->>}[ul]\ar @{>->}[dl]\ar[ll] _{f_{234}}
\\
     &&&  X_{0123} \ar@{->>}[ul]&& X_{1234}\ar @{>->}[dl] \ar@{->>}[ul]\ar[ll] _{f_{1234}}\\
      &&&&  \ar @/^4pc/@{->>}[uuuullll]_A
X_{01234}  \ar @{->>}[ul]
}
\end{gather*}

\vspace{-2mm}

We may now erase everything in the middle of the diagram to obtain the desired factorization
\begin{gather*}
\begin{split}
\xymatrix{
X_0  && X_1 \ar[ll]_{f_{1}}  &&  X_2
\ar[ll]_{f_{2}} &&    X_3 \ar[ll]_{f_{3}}
&&  X_4 \ar[ll]_{f_{4}}
\ar@/^4pc/@{>->}[ddddllll]_B
\\
 \\
\\
   \\
      &&&&  \ar @/^4pc/@{->>}[uuuullll]_A
X_{01234}
}
\end{split}\tag*{\qed}
\end{gather*}
\renewcommand{\qed}{}
\end{proof}

\begin{Remark}
Sylvain Cappell has pointed out the similarity of this result to ideas of J.H.C.~Whi\-tehead on simple homotopy theory, where maps are factored in to collapses and expansions.
And Thomas Kragh noted resemblances to the theory of Waldhausen categories; even the notation of decorated arrows is
similar.
\end{Remark}

\section{Rigid monoidal structures}
\label{sec-rigidmonoidal}

In this section, we introduce (rigid) monoidal structures and consider their compatibility with (highly) selective
structures.

\looseness=-1
Recall\footnote{We refer to~\cite{ba-ki:lectures} for a~detailed discussion of rigid monoidal categories.
An early reference for some of this material is~\cite[Chapter~1]{sa:categories}.} that a~\textit{monoidal structure}
on a~category ${\mathbf C}$ is a~functorial binary operation on objects and morphisms, which we will usually denote~by
$\otimes$.
This product operation is required to be associative and to admit a~unit object, usually denoted by $\mathbf{1}$; the
associativity and and unit identities (i.e.~equations) are assumed to hold up to specif\/ied natural isomorphisms which
satisfy certain further identities.
In some cases, these natural isomorphisms are simply identity morphisms, in which case the monoidal structure is known
as \textit{strict}.
Since the coherence theorem of MacLane~\cite{ma:categories} shows that every monoidal category is equivalent to a~strict
one, we will simplify our notation by assuming strictness; i.e., we will just identify a~product $A\otimes(B\otimes C)$
with $(A\otimes B)\otimes C$, and $A\otimes \mathbf{1}$ and $\mathbf{1} \otimes A$ with~$A$, rather than displaying the
identif\/ications explicitly as morphisms.

A monoidal structure is \textit{left rigid} if each object~$X$ is assigned a~\textit{left dual} object ${\overline{X}}$
equipped with \textit{unit} (sometimes called ``diagonal'') and \textit{counit} (sometimes called ``evaluation'')
morphisms $X\otimes {\overline{X}}\stackrel{\delta_X}{\longleftarrow}\mathbf{1}$ and
$\mathbf{1}\stackrel{\epsilon_X}{\longleftarrow}{\overline{X}}\otimes X$ for which each of the compositions below is
equal to an identity morphism
\begin{gather}
\label{eq-rigid1}
X\stackrel{1_X \otimes \epsilon_X}{\longleftarrow} X \otimes {\overline{X}} \otimes X \stackrel{\delta_X \otimes
1_X}{\longleftarrow X},
\\
\label{eq-rigid2}
{\overline{X}}\stackrel{\epsilon_X \otimes 1_{\overline{X}}}{\longleftarrow} {\overline{X}} \otimes X\otimes
{\overline{X}} \stackrel{1_{\overline{X}} \otimes \delta_X}{\longleftarrow {\overline{X}}}.
\end{gather}

There is a~similar def\/inition for right rigidity, and the structure is \textit{rigid} if it is both left and right
rigid.
In a~rigid monoidal category, there is a~natural bijection between ${\rm Hom}(X,Y\otimes Z)$ and ${\rm Hom}(X \otimes
{\overline{Y}}, Z)$ for any objects $X$, $Y$, and $Z$, in particular between ${\rm Hom}(X,Y)$ and ${\rm
Hom}(X\otimes{\overline{Y}},\mathbf{1})$.
Furthermore, the assignment of ${\overline{X}}$ to~$X$ extends to morphisms, def\/ining a~contravariant functor which is
an equivalence of categories.
The dual to a~morphism~$f$ will be denoted by $\overline f$.
The duality functor is also consistent with the monoidal structure: for any objects~$X$ and~$Y$, there is a~natural
isomorphism between ${\overline{Y}} \otimes {\overline{X}}$ and $\overline{X\otimes Y}$.
In addition, counits and units are related by the identity $\overline{\delta_X} = \epsilon_{\overline{X}}$.

\begin{Example}
\label{ex-relationsmonoidal}
The category $\mathbf{REL}$ has a~rigid monoidal structure in which the monoidal product of sets and relations is (and
will continue to be denoted as) the Cartesian product, and the unit object $\mathbf{1}$ has the empty set as its sole
element.
The (right and left) dual of any object is the object itself.
If we think of a~morphism $\xfy$ in $\mathbf{REL}$ as a~(multivalued, partially def\/ined, mapping) to~$X$ from~$Y$, the
corresponding element of ${\rm Hom}(X\times Y,\mathbf{1})$ is the subset of $X\times Y$ known as the graph of~$f$.

Similar def\/initions apply to the categories $\mathbf{MREL}$ and $\mathbf{SREL}$ of smooth and canonical relations.
\end{Example}

Guided by the example above, we will use the term \textit{graph of}~$f$ and the notation $\gamma_f$ to denote, for any
rigid monoidal category, the morphism to $X \otimes {\overline{Y}}$ from $ \mathbf{1}$ corresponding to $\xfy$.
In particular, the unit morphism $X\otimes {\overline{X}}\stackrel{\delta_X}{\leftarrow} \mathbf{1}$ is the graph of the
identity morphism of~$X$.

More explicitly, the graph $\gamma_f$ of $\xfy$ is the composition
\begin{gather*}
X\otimes {\overline{Y}} \stackrel{f\otimes 1_{{\overline{Y}}}}{\longleftarrow}Y\otimes {\overline{Y}}
\stackrel{\delta_Y}{\longleftarrow} \mathbf{1}.
\end{gather*}

On the other hand, given any~$\gamma$ to $X\otimes {\overline{Y}}$ from $\mathbf{1}$, we may def\/ine the morphism
$f_\gamma$ to~$X$ from~$Y$ as the composition
\begin{gather*}
X \stackrel {1_X\otimes \epsilon_Y}{\longleftarrow} X\otimes {\overline{Y}} \otimes Y \stackrel {\gamma \otimes
1_Y}{\longleftarrow} Y.
\end{gather*}

Showing that these operations between morphisms and their graphs are inverse to one another is an exercise in applying
the def\/ining identities of a~monoidal category.
Since the solution to this exercise is omitted from all the references we have found, we include it here for
completeness.

\begin{Proposition}
With the definitions above, $f_{\gamma_f}=f$ for any $\xfy$ and $\gamma_{f_\gamma}= \gamma$ for any $X\otimes
{\overline{Y}} \stackrel{\gamma}{\leftarrow}\mathbf{1}$.
\end{Proposition}
\begin{proof}
By the def\/initions,
\begin{gather*}
f_{\gamma_f}=(1_X\otimes \epsilon_Y) (\gamma_f \otimes 1_Y),
\end{gather*}
which equals
\begin{gather*}
(1_X\otimes \epsilon_Y)(f\otimes 1_{{\overline{Y}}}\otimes 1_Y)(\delta_Y\otimes 1_Y).
\end{gather*}
Carrying out the compositions with identity operators, we may combine the f\/irst two factors, giving us
\begin{gather*}
(f\otimes \epsilon_Y)(\delta_Y\otimes 1_Y).
\end{gather*}
Now we rewrite~$f$ as $f 1_Y$, giving
\begin{gather*}
f (1_Y \otimes \epsilon_Y)(\delta_Y \otimes 1_Y).
\end{gather*}
The last two factors collapse to $1_Y$ by one of the rigidity identities, leaving us with~$f$.

For the other direction, we have
\begin{gather*}
\gamma_{f_\gamma} =(f_\gamma \otimes 1_{\overline{Y}})\delta_Y,
\end{gather*}
which in terms of~$\gamma$ itself is
\begin{gather*}
(1_X\otimes \epsilon_Y\otimes1_{\overline{Y}})(\gamma\otimes 1_Y\otimes 1_{\overline{Y}})\delta_Y.
\end{gather*}
Combining the last two factors gives
\begin{gather*}
(1_X\otimes \epsilon_Y\otimes1_{\overline{Y}})(\gamma\otimes\delta_Y),
\end{gather*}
which can then be expanded as
\begin{gather*}
(1_X\otimes \epsilon_Y\otimes1_{\overline{Y}})(1_X\otimes 1_{\overline{Y}}\otimes\delta_Y)\gamma.
\end{gather*}
Using the rigidity identity $(\epsilon_Y\otimes1_{\overline{Y}})(1_{\overline{Y}}\otimes\delta_Y)=1_{\overline{Y}}$, we
arrive at~$\gamma$.
\end{proof}

Now we will look at the composition of morphisms in terms of their graphs.

\begin{Definition}
Given morphisms $\gamma_1$ and $\gamma_2$ to $X\otimes {\overline{Y}}$ and $Y \otimes {\overline{Z}}$ respectively from
$\mathbf{1}$, their \textit{reduced product} is the morphism $\gamma_1 \ast \gamma_2$ to $X\otimes {\overline{Z}}$ from
$\mathbf{1}$ given~by
\begin{gather*}
\gamma_1 \ast \gamma_2 = (1_X \otimes \epsilon_Y\otimes 1_{\overline{Z}}) (\gamma_1\otimes \gamma_2).
\end{gather*}
\end{Definition}
The following result must be well-known, but we could not f\/ind it in the literature.

\begin{Proposition}
For any morphisms $\xfy$ and $\ygz$ in a~rigid monoidal category, the graph of their composition is equal to the reduced
product $\gamma_f \ast \gamma_g$ of their graphs.
\end{Proposition}
\begin{proof}
It suf\/f\/ices to show that the morphism to~$X$ from~$Y$ associated with $\gamma_f \ast \gamma_g$ is the composition $fg$.

This associated morphism is
\begin{gather*}
(1_X \otimes \epsilon_Z)((\gamma_f\ast\gamma_g)\otimes 1_Z).
\end{gather*}
By the def\/inition of the reduced product and the fact that $1_Z= 1_Z 1_Z$, we get
\begin{gather*}
(1_X \otimes \epsilon_Z) (1_X \otimes \epsilon_Y\otimes 1_{\overline{Z}}\otimes 1_Z) (\gamma_f\otimes \gamma_g\otimes
1_Z).
\end{gather*}
The composition of the f\/irst two factors goes through $X\otimes{\overline{Z}}\otimes Z$.
We can go through $X\otimes {\overline{Y}}\otimes Y$ instead to get
\begin{gather*}
(1_X \otimes \epsilon_Y)(1_X\otimes 1_{\overline{Y}}\otimes 1_Y \otimes \epsilon_Z)(\gamma_f\otimes \gamma_g\otimes
1_Z).
\end{gather*}

We can factor the third factor through $Y\otimes {\overline{Z}} \otimes Z$ to get
\begin{gather*}
(1_X \otimes \epsilon_Y)(1_X\otimes 1_{\overline{Y}}\otimes 1_Y \otimes \epsilon_Z) (\gamma_f \otimes 1_Y \otimes
1_Z\otimes 1_{\overline{Z}})(\gamma_g \otimes 1_Z).
\end{gather*}
Monoidal identities allow us to simplify the product of the middle two factors to get
\begin{gather*}
(1_X \otimes \epsilon_Y)(\gamma_f\otimes 1_Y)(1_Y\otimes \epsilon_Z)(\gamma_g \otimes 1_Z).
\end{gather*}
The product of the f\/irst two factors is now~$f$, and that of the last two is~$g$, so we are done.
\end{proof}

\begin{Remark}
The reader may f\/ind it useful to f\/ill in the sources and targets of the morphisms in the argument above, perhaps drawing
a~diagram.
On the other hand, suppressing the identity morphisms produces the following nice condensed version of the proof
\begin{gather*}
\epsilon_Z(\gamma_f \ast \gamma_g) = \epsilon_Z \epsilon_Y (\gamma_f \otimes \gamma_g) = \epsilon_Y \epsilon_Z (\gamma_f
\otimes \gamma_g) =\epsilon_Y \gamma_f \epsilon_Z \gamma_g = fg.
\end{gather*}
\end{Remark}

In any rigid monoidal category, the endomorphisms of the unit object play a~special role.
First of all, they form a~commutative monoid under either of two operations, which turn out to coincide\footnote{The
following proof is sometimes known as the Eckmann--Hilton argument; see \url{http://math.ucr.edu/home/baez/week258.html}.}.

\begin{Proposition}
\label{prop-endounitob}
If~$f$ and~$g$ are endomorphisms of $\mathbf{1}$ in a~rigid monoidal category, then $f\otimes g$, $fg$, $g\otimes f$,
and $gf$ are all equal.
Thus, ${\rm Hom} (\mathbf{1},\mathbf{1})$ forms a~commutative monoid under either operation, with $1_\mathbf{1}$ as its
identity element.
\end{Proposition}
\begin{proof}
It is clear that $1_\mathbf{1}$ is an identity element for either operation.
Now
\begin{gather*}
f\otimes g = 1_\mathbf{1} f \otimes g 1_\mathbf{1} = (1_\mathbf{1} \otimes g)(f\otimes 1_\mathbf{1})= gf.
\end{gather*}
But we also have
\begin{gather*}
f\otimes g = f 1_\mathbf{1} \otimes 1_\mathbf{1} g = fg.
\end{gather*}
Interchanging~$f$ and~$g$ gives the f\/inal equality, $g\otimes f = gf$.

Another way of seeing the equality of composition and monoidal product in ${\rm Hom} (\mathbf{1},\mathbf{1})$ is to
observe that each such endomorphism is {\em equal} to its graph.
\end{proof}

An argument similar to the proof of Proposition~\ref{prop-endounitob}, using strictness of the monoidal structure to
identify endomorphisms of $\mathbf{1}\otimes X$ with those of $X\otimes \mathbf{1}$, shows:

\begin{Proposition}
The monoid of endomorphisms of the unit object in a~monoidal category ${\mathbf C}$ acts on ${\mathbf C}$ by monoidal
product from the left or right.
These actions are equal if the monoidal product is symmetric.
\end{Proposition}

There is also a~useful \textit{trace} morphism~$\tau$ to ${\rm Hom}(\mathbf{1},\mathbf{1})$ from ${\rm Hom}(X,X)$ for
any~$X$.
It is def\/ined by $\tau(f) = \epsilon_X \gamma_f$.

We turn now to the interaction of (rigid) monoidal structures and (highly) selective structures.

\begin{Definition}
\label{dfn-selectivemonoidal}
A~\textit{selective monoidal category} is a~monoidal category together with a~selective structure such that the classes
of suave morphisms and of congenial pairs are closed under the monoidal product operation.

If, in addition, the category ${\mathbf C}$ is highly selective, and the subcategories of reductions and coreductions
are closed under the monoidal product, then we call ${\mathbf C}$ a~\textit{highly selective monoidal category}.

A \textit{selective rigid monoidal category} is a~selective monoidal category equipped with a~rigid structure such that:
\begin{enumerate}
\itemsep=0pt
\item The unit and counit morphisms are suave.
\item The duality functor is selective.
\item If $X\stackrel{f}{\leftarrow} Y\otimes Z$ is suave, then the compositions $(f\otimes 1_{\overline{Z}})(1_Y\otimes
\delta_Z)$ and $(1_{\overline{Y}}\otimes f) (\delta_{\overline{Y}}\otimes 1_Z)$ are both congenial.
\end{enumerate}
If, in addition, the category is highly selective, it is a~\textit{highly selective rigid monoidal category} if the unit
and counit morphisms are coreductions and reductions respectively, and the duality functor exchanges reductions and
coreductions.
\end{Definition}

The congeniality of other compositions follows immediately by duality from the axioms above.
We state this as a~proposition and omit the simple proof.
\begin{Proposition}

\label{prop-congenialepsilon}
In a~selective rigid monoidal category, if $X\times Y \stackrel{g}{\leftarrow}Z$ is suave, then the compositions
$(\epsilon_X\otimes 1_Y)(1_{\overline{X}}\otimes g)$ and $(1_X\otimes\epsilon_Y)(g\otimes 1_{\overline{Y}})$ are
congenial.

In particular, the compositions~\eqref{eq-rigid1} and~\eqref{eq-rigid2},
\begin{gather*}
X\stackrel{1_X \otimes \epsilon_X}{\longleftarrow} X \otimes {\overline{X}} \otimes X \stackrel{\delta_X \otimes
1_X}{\longleftarrow X},
\qquad
{\overline{X}}\stackrel{\epsilon_X \otimes 1_{\overline{X}}}{\longleftarrow} {\overline{X}} \otimes X\otimes
{\overline{X}} \stackrel{1_{\overline{X}} \otimes \delta_X}{\longleftarrow {\overline{X}}},
\end{gather*}
required to give identity morphisms in the def\/inition of a~rigid monoidal category are congenial.
\end{Proposition}

\begin{Example}
\label{ex-relationsselectiverigidity}
The category $\mathbf{REL}$ of relations, with the structures described in
Examp\-les~\ref{ex-selectivesuave},~\ref{ex-relhighlyselective}, and~\ref{ex-relationsmonoidal}, is a~highly selective
rigid monoidal category.
Conditions 1 and 2 in Def\/inition~\ref{dfn-selectivemonoidal} are obviously satsif\/ied.
For condition 3, given~$f$, one forms the composition $(f\otimes 1_{\overline{Z}})(1_Y\otimes \delta_Z)$ by taking all
9-tuples of the form
\begin{gather*}
((x,z),(y',z',z),(y,z'',z''),y)\in (X\times {\overline{Z}}) \times (Y\times Z \times {\overline{Z}}) \times (Y\times
Z\times{\overline{Z}}) \times Y
\end{gather*}
for which $(x,(y,z))\in f$ and $(y',z',z)=(y,z'',z'')$, and then projecting to $((x,z),y) \in (X\times
{\overline{Z}})\times Y $.
The equalities $y'=y$ and $z''=z'=z$ imply that the intermediate elements $y'$, $z'$, and~$z$ are all determined
by~$x$,~$z$, and~$y$, so the composition is monic, i.e.~congenial.

Congeniality of the other composition in condition 3 is verif\/ied in the same way.
\end{Example}

Before developing further consequences of selectivity combined with rigidity, we make a~short digression to discuss
string diagrams.
These diagrams provide a~graphical notation for morphisms in rigid monoidal categories which facilitates the proofs of
many identities in such categories.
For simplicity, we will restrict attention to symmetric monoidal categories\footnote{String diagrams are also useful in
the absence of the symmetry assumption, but care must be taken to distinguish between left and right duals.}.

\looseness=1
One depicts objects (either the targets or sources of morphisms) by directed vertical line segments: an upwards directed
line labelled by~$X$ denotes~$X$; dually, a~downwards directed line labelled by~$X$ denotes the dual, ${\overline{X}}$.
One depicts a~morphism $X_1\otimes\dots \otimes X_n\xleftarrow{f}Y_1\otimes\dots\otimes Y_m$ as\vspace*{-3mm}
\begin{gather}
\label{eq:StringMorph}
\begin{split}
&\begin{tikzpicture}
\node (x1) at (-.5,1.5){};
\node (ex) at (0,.75) {$\cdots$};
\node (xn) at (.5,1.5){};
\node (f) at (0,0) [shape=rectangle,draw] {$\quad f\quad $};
\node (y1) at (-.5,-1.5) {};
\node (ey) at (0,-.75) {$\cdots$};
\node (ym) at (.5,-1.5) {};
\draw[->] ($(f.north)+(-.5,0)$) -- node {$X_1$} (x1);
\draw[->] ($(f.north)+(.5,0)$) -- node [swap]{$X_n$} (xn);
\draw[->] (y1) -- node {$Y_1$} ($(f.south)+(-.5,0)$);
\draw[->] (ym) -- node [swap] {$Y_m$} ($(f.south)+(.5,0)$);
\node (or) at (2,0) {$\dot=$};
\node (x1) at (3.5,1.5){};
\node (ex) at (4,.75) {$\cdots$};
\node (xn) at (4.5,1.5){};
\node (f) at (4,0) [shape=rectangle,draw] {$\quad f\quad $};
\node (y1) at (3.5,-1.5) {};
\node (ey) at (4,-.75) {$\cdots$};
\node (ym) at (4.5,-1.5) {};
\draw[->] ($(f.north)+(-.5,0)$) -- node {$X_1$} (x1);
\draw[->] ($(f.north)+(.5,0)$) -- node [swap]{$X_n$} (xn);
\draw[<-] (y1) -- node {$\overline{Y}_1$} ($(f.south)+(-.5,0)$);
\draw[<-] (ym) -- node [swap] {$\overline{Y}_m$} ($(f.south)+(.5,0)$);
\end{tikzpicture}
\end{split}
\end{gather}

\vspace*{-3mm}

\noindent
Here, we use the notation $\dot =$ to signify that two diagrams represent the same morphism in the underlying category.

To depict the monoidal product of objects/morphisms, one places those objects/morphisms horizontally side by side.
For composition of morphisms, the morphisms are stacked vertically, and the segments denoting the source/target objects
are connected in the appropriate way.
For example, given morphisms $X_1\xleftarrow{f} Y_1$, $X_2\xleftarrow{g}Y_2$, and $Y_1\otimes Y_2\xleftarrow{h} Z$, the
morphism $(f\otimes g)h$ is depicted as\vspace*{-3mm}
\begin{gather*}
\begin{tikzpicture}
\node (x1) at (-.5,1.5) {};
\node (x2) at (.5,1.5) {};
\node (fgh) at (0,0) [shape=rectangle,draw] {$(f\otimes g)h$};
\node (z) at (0,-1.5) {};
\draw[->] ($(fgh.north)+(-.5,0)$) -- node {$X_1$} (x1);
\draw[->] ($(fgh.north)+(.5,0)$) -- node [swap]{$X_2$} (x2);
\draw[->] (z) -- node {$Z$} (fgh);
\node (x1') at (4,2.25) {};
\node (x2') at (5,2.25) {};
\node (f) at (4,.75) [shape=rectangle,draw] {$f$};
\node (g) at (5,.75) [shape=rectangle,draw] {$g$};
\node (h) at (4.5,-.75) [shape=rectangle,draw] {$\quad h\quad$};
\node (z) at (4.5,-2.25) {};
\draw[->] (f) -- node {$X_1$} (x1');
\draw[->] (g) -- node [swap] {$X_2$} (x2');
\draw[->] ($(h.north)+(-.5,0)$) -- node {$Y_1$} (f);
\draw[->] ($(h.north)+(.5,0)$) -- node[swap] {$Y_2$} (g);
\draw[->] (z) -- node {$Z$} (h);
\node (eq) at (2,0) {$\dot =$};
\end{tikzpicture}
\end{gather*}
For the identity morphisms, one omits the boxes:
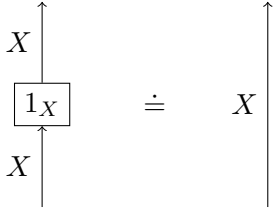
\vspace*{-3mm}
\begin{gather*}
\begin{tikzpicture}
\node (x) at (0,1.5){};
\node (f) at (0,0) [shape=rectangle,draw] {$1_X$};
\node (y) at (0,-1.5) {};
\draw[->] (f) -- node {$X$} (x);
\draw[->] (y) -- node {$X$} (f.south);
\node (eq) at (1.5,0) {$\dot =$};
\node (x) at (3,1.5){};
\node (y) at (3,-1.5) {};
\draw[->] (y) -- node {$X$} (x);
\end{tikzpicture}
\end{gather*}
One also omits the boxes for the unit and counit morphisms, as well as for the segment corresponding to the terminal
object
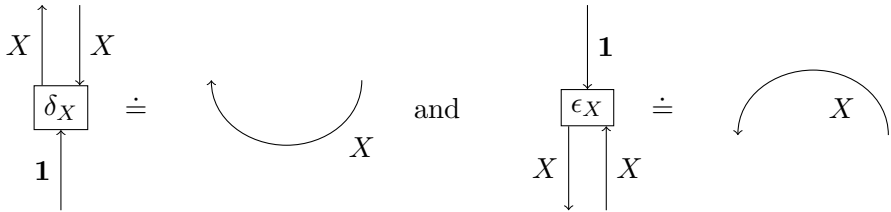
\vspace*{-3mm}
\begin{gather*}
\begin{tikzpicture}
\node (x) at (-.25,1.5) {};
\node (bx) at (.25,1.5) {};
\node (d) at (0,0) [shape=rectangle,draw] {$\delta_X$};
\node (1) at (0,-1.5) {};
\draw[->] ($(d.north)+(-.25,0)$) -- node {$X$} (x);
\draw[<-] ($(d.north)+(.25,0)$) -- node [swap]{$X$} (bx);
\draw[->] (1) -- node {$\mathbf{1}$} (d);
\node (eq) at (1,0) {$\dot=$};
\node (x') at (2,.5) {};
\node (bx') at (4,.5) {};
\draw[->] (bx'.south) to [out=270,in=0] node {$X$} (3,-.5) to [out=180,in=270] (x'.south);
\node (and) at (5,0) {and};
\node (x) at (6.75,-1.5) {};
\node (bx) at (7.25,-1.5) {};
\node (d) at (7,0) [shape=rectangle,draw] {$\epsilon_X$};
\node (1) at (7,1.5) {};
\draw[->] ($(d.south)+(-.25,0)$) -- node [swap]{$X$} (x);
\draw[<-] ($(d.south)+(.25,0)$) -- node {$X$} (bx);
\draw[->] (1) -- node {$\mathbf{1}$} (d);
\node (eq) at (8,0) {$\dot=$};
\node (x') at (9,-.5) {};
\node (bx') at (11,-.5) {};
\draw[->] (bx'.north) to [out=90,in=0] node {$X$} (10,.5) to [out=180,in=90] (x'.north);
\end{tikzpicture}
\end{gather*}
With this notation, one can write the rigidity axioms~\eqref{eq-rigid1} and~\eqref{eq-rigid2} as
\begin{gather*}
\begin{tikzpicture}
\draw[->] (1,-.75) -- (1,0) to [out=90,in=0] (.5,.5) to [out=180,in=90] (0,0) to [out=270,in=0] (-.5,-.5) to [out=180,in=270] (-1,0) -- node {$X$} (-1,.75);
\node at (2,0) {$\dot =$};
\draw[->] (3,-.75) -- node {$X$} (3,.75);
\node at (4,0) {$\dot =$};
\draw[->] (5,-.75) -- (5,0) to [out=90,in=180] (5.5,.5) to [out=0,in=90] (6,0) to [out=270,in=180] (6.5,-.5) to [out=0,in=270] (7,0) -- node {$X$} (7,.75);
\end{tikzpicture}
\end{gather*}

In a~string diagram~$D$ representing some composite morphism $f_D$, we will refer to the paths composed of identity
morphisms, units, and counits as {\em strings}.

To signify that the composition of $X\xleftarrow{f}Y$ and $Y\xleftarrow{g}Z$ is congenial, we mark the string diagram
with a~$\bullet$:
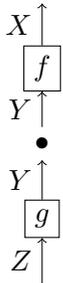
\vspace*{-3mm}
\begin{gather*}
\begin{tikzpicture}
\node (x) at (0,1){};
\node (f) at (0,0) [shape=rectangle,draw] {$f$};
\node (y) at (0,-1) {$\bullet$};
\node (g) at (0,-2) [shape=rectangle,draw] {$g$};
\node (z) at (0,-3){};
\draw[->] (f) -- node {$X$} (x);
\draw[->] (y) -- node {$Y$} (f);
\draw[->] (g) -- node {$Y$} (y);
\draw[->] (z) -- node {$Z$} (g);
\end{tikzpicture}
\end{gather*}

Then Axiom~2 from Def\/inition~\ref{dfn-selective} can be written as:

$2'$.
{\em Whenever $X\xleftarrow{g} Y$ is suave, the following compositions are congenial:}\vspace*{-3mm}
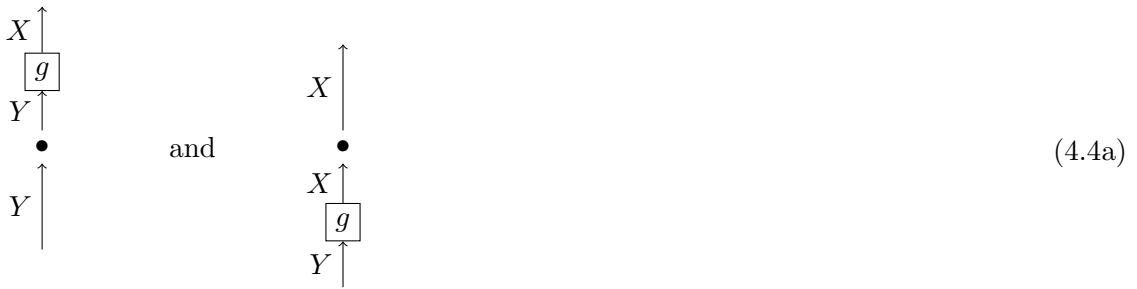
\begin{subequations}
\begin{gather}
\label{eq:StrIdCong}
\begin{split}
& \begin{tikzpicture}
\node (x) at (0,1){};
\node (f) at (0,0) [shape=rectangle,draw] {$g$};
\node (y) at (0,-1) {$\bullet$};
\node (z) at (0,-2.5){};
\draw[->] (f) -- node {$X$} (x);
\draw[->] (y) -- node {$Y$} (f);
\draw[->] (z) -- node {$Y$} (y);
\node at (2,-1) {and};
\node (x) at (4,.5){};
\node (y) at (4,-1) {$\bullet$};
\node (g) at (4,-2) [shape=rectangle,draw] {$g$};
\node (z) at (4,-3){};
\draw[->] (y) -- node {$X$} (x);
\draw[->] (g) -- node {$X$} (y);
\draw[->] (z) -- node {$Y$} (g);
\end{tikzpicture}
\end{split}
\end{gather}

Axiom 3 from Def\/inition~\ref{dfn-selectivemonoidal} can be written as:

$3'$.
{\em Whenever $X\xleftarrow{f}Y\otimes Z$ is suave, the following compositions are congenial:}\vspace*{-3mm}
\begin{gather}
\label{eq:StrUnCong}
\begin{split}
& \begin{tikzpicture}
\node (x) at (0,1.25){};
\node (tby) at (-1.25,1.25){};
\node (f) at (0,0) [shape=rectangle,draw] {$\; f \;$};
\node (y) at (-.25,-1.25) {$\bullet$};
\node (z) at (.25,-1.25) {$\bullet$};
\node (by) at (-1.25,-1.25) {$\bullet$};
\node (z') at (.25,-2.25) {};
\draw[->] (f) -- node {$X$} (x);
\draw[->] (y) -- node {$Y$} ($(f.south)+(-.25,0)$);
\draw[->] (z) -- node [swap]{$Z$} ($(f.south)+(.25,0)$);
\draw[->] (tby) -- node [swap]{$Y$} (by);
\draw[->] (by.south) to [out=270,in=180] node [swap]{$Y$} (-.75,-2) to [out=0,in=270] (y.south);
\draw[->] (z') -- node [swap]{$Z$} (z);
\node at (2,0) {and};
\node (x) at (4,1.25){};
\node (tbz) at (5.25,1.25){};
\node (f) at (4,0) [shape=rectangle,draw] {$\; f \;$};
\node (y) at (3.75,-1.25) {$\bullet$};
\node (z) at (4.25,-1.25) {$\bullet$};
\node (bz) at (5.25,-1.25) {$\bullet$};
\node (y') at (3.75,-2.25) {};
\draw[->] (f) -- node {$X$} (x);
\draw[->] (y) -- node {$Y$} ($(f.south)+(-.25,0)$);
\draw[->] (z) -- node [swap]{$Z$} ($(f.south)+(.25,0)$);
\draw[->] (tbz) -- node {$Z$} (bz);
\draw[->] (bz.south) to [out=270,in=0] node {$Z$} (4.75,-2) to [out=180,in=270] (z.south);
\draw[->] (y') -- node {$Y$} (y);
\end{tikzpicture}
\end{split}
\end{gather}
\end{subequations}

\looseness=-1
Suppose that $(D_1,D_2,\dots,D_n)$ is a~composable~$n$-tuple of string diagrams representing a~composable~$n$-tuple of
composite morphisms $(f_{D_1},f_{D_2},\dots,f_{D_n})$.
Then the two axioms above imply that $(f_{D_1},\dots,f_{D_n})$ is a~congenial~$n$-tuple if the boxes in each diagram
$D_i$ are all labelled by suave morphisms, and if one may redirect\footnote{Note that, if one changes the direction of
a~string labelled by the object~$X$, one must also change the label to its dual, ${\overline{X}}$.} some of the strings
in each diagram $D_i$ in such a~way that
\begin{enumerate}\itemsep=0pt
\item[s1)] no string enters a~box\footnote{In~\eqref{eq:StringMorph}, the right hand side is an example of a~string
diagram in which no string enters a~box, while the left hand side is not.}, and
\item[s2)] all strings compose tip to tail.
\end{enumerate}

For example, the diagrams in~\eqref{eq:StrIdCong} and~\eqref{eq:StrUnCong} do not satisfy conditions (s1) and (s2).
However, after redirecting the strings labelled by the objects~$Y$ and~$Z$ (and relabelling them ${\overline{Y}}$ and
${\overline{Z}}$) the resulting diagrams do satisfy~(s1) and~(s2).

We return now to developing some properties of rigidity combined with selectivity.
The reader is invited to represent some of our proofs by string diagrams.

\begin{Proposition}
\label{prop-suavegraph}
In a~selective rigid monoidal category, a~morphism is suave if and only if its graph is.
\end{Proposition}

\begin{proof}
The graph $\gamma_f$ of $\xfy$ is the composition $(f\otimes 1_{\overline{Y}})\delta_Y$.
This can also be written as $(f\otimes 1_{\overline{Y}})(1_\mathbf{1}\otimes\delta_Y)$.
By Def\/inition~\ref{dfn-selectivemonoidal}, this composition is congenial if~$f$ is suave, so $\gamma_f$ is suave.

On the other hand, the morphism $f_\gamma$ with graph~$\gamma$ is
\begin{gather*}
X \stackrel {1_X\otimes \epsilon_Y}{\longleftarrow} X\otimes {\overline{Y}} \otimes Y \stackrel {\gamma \otimes
1_Y}{\longleftarrow} Y.
\end{gather*}
If~$\gamma$ is suave, the composition is congenial by Proposition~\ref{prop-congenialepsilon}, so $f_\gamma$ is suave.
\end{proof}

This leads to a~way of passing from rigidity and selectivity to high selectivity.

\begin{Proposition}
\label{prop-factorization}
Let be ${\mathbf C}$ be a~selective rigid monoidal category with distinguished rigid monoidal subcategories of
reductions and coreductions satisfying the first two conditions in De\-fi\-nition~{\rm \ref{dfn-highly}}.
Suppose further that morphisms from and to unit objects are coreductions and reductions respectively.
Then the factorization axiom is also satisfied, and ${\mathbf C}$ is a~highly selective rigid monoidal category.
\end{Proposition}

\begin{proof}
The last equation in the proof of Proposition~\ref{prop-suavegraph} provides the required factori\-za\-tion.
\looseness=2
\end{proof}

We turn now to the WW construction in the (rigid) monoidal setting.

If ${\mathbf C}$ is any selective monoidal category, we can try to put a~monoidal structure on ${\rm WW}({\mathbf C})$ as
follows.
The monoidal product on objects and the unit objects are the same as in ${\mathbf C}$.
To def\/ine the monoidal product of $\rm WW$-morphisms, we choose representative sequences for them and form their their
monoidal product entry by entry to get a~representative of the product.
(This is one place where it is useful to represent the morphisms by inf\/inite sequences.)

To show that this def\/inition is correct, we will verify the monoidal product of $\rm WW$-morphisms is well-def\/ined,
associative, and functorial.

\begin{Proposition}
The product operation described above is well defined, i.e., the monoidal product of two morphisms in ${\rm WW}({\mathbf C})$
is independent of the choice of representatives.
\end{Proposition}

\begin{proof}
Let $(f_i)$ and $(f'_i)$ be sequences representing morphisms in ${\rm WW}({\mathbf C})$.
We will show that the insertion of a~unit morphism or the collapsing of a~congenial pair in either sequence does not
change the equivalence class of the entry-by-entry monoidal product.

We deal f\/irst with the insertion of identity morphisms.
For simplicity of notation, we will omit their subscripts, which are determined by the context.
Suppose now that the appropriate identity morphism is inserted after some $f_i$ and everything else is shifted one step
to the right.
Then, in the monoidal product sequence, all entries to the left of $f_i\otimes f'_i$ remain the same, the $i+1$'st entry
becomes $1\otimes f'_{i+1}$, the next are $f_{i+1}\otimes f'_{i+2}$, $f_{i+2}\otimes f'_{i+3}$, etc., eventually
becoming the identity morphism of the source object.
To convert the original sequence to this form, we f\/irst factor each $f_{i+k}\otimes f'_{i+k+1}$, replacing it by the
congenial pair $(1\otimes f'_{i+k+1},f_{i+k}\otimes 1)$, until we arrive at the identity morphisms.
Leaving the entry $(1\otimes f'_{i+1})$ in place and composing pairwise the remaining ones on the right yields the
desired sequence.
An analogous argument applies if we shift to the left after insertion, or if we insert the identity morphism in the
second monoidal factor.

Since we have dealt with the insertion of identity morphisms, to handle the composition of congenial pairs, it suf\/f\/ices
to consider the situation where some congenial $(f_i,f_{i+1})$ is replaced by $(f_i f_{i+1},1_{X_{i+1}})$.

It suf\/f\/ices to look at the segment $(f_i \otimes f'_i,f_{i+1}\otimes f'_{i+1})$ of the monoidal product sequence, since
the rest is unchanged.
Using the facts that monoidal products of congenial sequences are congenial and that composable sequences with a~unit
morphism as entry are congenial, we may replace the segment $(f_i \otimes f'_i,f_{i+1}\otimes f'_{i+1})$ by the
equivalent segment $(f_i \otimes f'_i, f_{i+1} \otimes 1, 1 \otimes f'_{i+1})$, and then by $(f_i f_{i+1} \otimes f'_i,
1 \otimes f'_{i+1})$.

Of course, this argument applies as well if the composition is done on the second monoidal factor, and our proof is
complete.
\end{proof}

\begin{Proposition}
The monoidal product defined above is associative.
It is also a~bifunctor: if $[f_j]\in {\rm WW}({\mathbf C})(X,Y)$, $[g_j] \in {\rm WW}({\mathbf C})(Y,Z)$, $[f'_j]\in {\rm WW}({\mathbf
C})(X',Y')$, and $[g'_j] \in {\rm WW}({\mathbf C})(Y',Z')$, then $([f_j]\otimes [f'_j])\circ ([g_j]\otimes [g'_j]) =[f_j]\circ
[g_j]\otimes [f'_j]\circ [g'_j]$.
\end{Proposition}
\begin{proof}
Associativity follows immediatively from associativity of the monoidal product in ${\mathbf C}$.
For functoriality, we use the idea in Remark~\ref{rmk-shifting}, choosing representative sequences and then shifting
them so that the supports of the~$f$ and~$g$ sequences are contained in the negative and postive integers respectively.
The same is then true for the sequences $(f_j\otimes f'_j)$ and $(g_j \otimes g'_j)$, from which the desired result
follows.
\end{proof}

If ${\mathbf C}$ is highly selective monoidal, then the reductions and coreductions are monoidal sub\-ca\-te\-gories of
${\rm WW}({\mathbf C})$.

As for rigidity, suppose that ${\mathbf C}$ is a~selective rigid monoidal category.
We have seen that ${\rm WW}({\mathbf C})$ is again monoidal.
To extend the rigid structure, we simply let the duality operation on objects be that of ${\mathbf C}$, with the
morphisms $\delta_X$ and $\epsilon_X$ being their images under the embedding~$c$ from suave morphisms to WW-morphisms.
Since we are assuming that the compositions in the rigidity axioms are congenial, the axioms continue to hold in
${\rm WW}({\mathbf C})$.
We may therefore conclude:

\begin{Theorem}
If ${\mathbf C}$ is a~selective $($rigid$)$ monoidal category, then ${\rm WW}({\mathbf C})$ is a~$($rigid$)$ monoidal category.
If ${\mathbf C}$ is, in addition, highly selective, then ${\rm WW}({\mathbf C})$ contains distinguished monoidal subcategories
of reductions and coreductions.
\end{Theorem}

\section{Hypergraphs}
If $\gamma_f$ and $\gamma_g$ are suave graphs in a~selective rigid monoidal category ${\mathbf C}$,~$f$ and~$g$ are
suave morphisms, but the composition $fg$ is generally not suave.
This means that the reduced product $\gamma_f \ast \gamma_g$ is not suave either.
To describe the reduced product in terms of suave morphisms, it is natural to pass to ${\rm WW}({\mathbf C})$, with its
inherited rigid monoidal structure.
In particular, we still have in ${\rm WW}({\mathbf C})$ a~bijection between morphisms and their graphs.

We will assume, further, that ${\mathbf C}$ is highly selective rigid monoidal.
By Theorem~\ref{thm-two morphisms}, any morphism $X\otimes{\overline{Y}} \leftarrow\mathbf{1}$is represented~by
a~composition of the form
\begin{gather*}
X\otimes {\overline{Y}} \twoheadleftarrow R \leftarrowtail \mathbf{1}
\end{gather*}
in the original category ${\mathbf C}$.
We will call such a~diagram a~\textit{hypergraph}.
Note that dif\/ferent hypergraphs may be $\rm WW$-equivalent, so that this representation of a~$\rm WW$-morphism is not unique.
For example, the identity on~$X$ is represented by the hypergraph $X\otimes {\overline{X}} \twoheadleftarrow X\otimes
{\overline{X}} \stackrel{\delta_X}{\leftarrowtail} \mathbf{1}$, but composing this with itself according to the rule
below gives a~dif\/ferent representative.

To compose hypergraphs $X\otimes{\overline{Y}} \twoheadleftarrow R \leftarrowtail \mathbf{1}$ and $Y\otimes
{\overline{Z}} \twoheadleftarrow S \leftarrowtail \mathbf{1}$, we start with the four-step sequence
\begin{gather*}
X\otimes {\overline{Z}} \xleftarrow{1_X \otimes \epsilon_Y\otimes 1_Z} X\otimes {\overline{Y}}\otimes Y\otimes
{\overline{Z}} \twoheadleftarrow R \otimes S \leftarrowtail \mathbf{1}\otimes\mathbf{1}\leftarrowtail \mathbf{1}.
\end{gather*}
The f\/irst two and last two arrows form congenial pairs, so we may replace them by their compositions to obtain the
hypergraph
\begin{gather*}
X\otimes {\overline{Z}} \twoheadleftarrow R\otimes S \leftarrowtail \mathbf{1},
\end{gather*}
which represents the reduced product of the original two.
Since ${\rm WW}({\mathbf C})$ is rigid monoidal, the WW-equivalence class of this reduced product depends only on that of the 
hypergraphs being combined.

Another useful property of hypergraphs is just a~special case of Proposition~\ref{prop-compwithisomo}.
\begin{Proposition}
\label{prop-hyperisomo}
Hypergraphs $X\otimes{\overline{Y}} \stackrel{C}{\twoheadleftarrow} R \stackrel{L}{\leftarrowtail} \mathbf{1}$ and
$X\otimes{\overline{Y}} \twoheadleftarrow R' \leftarrowtail \mathbf{1}$ are $\rm WW$-equivalent if they are isomorphic in
the sense that there is a~suave isomorphism~$\phi$ to~$R$ from $R'$ such that $\phi L' = L$ and $C' \phi = C$.
\end{Proposition}

\begin{Remark}
There is a~2-category with the hypergraphs as morphisms and the isomorphisms of hypergraphs (or more general morphisms,
where~$\phi$ is not invertible, but the relevant compositions are congenial) as 2-morphisms.
But we will see many examples of $\rm WW$-equivalent hypergraphs which are not isomorphic.
It would be very interesting to include these more general equivalences as $2$-morphisms in a~larger $2$-category.
Something like this is done for a~special case in~\cite{we-wo:functoriality}.
\end{Remark}

An important class of hypergraphs consists of those for which the source is the unit object.
In this case, a~hypergraph is simply a~representative $X\twoheadleftarrow R \leftarrowtail \mathbf{1}$ of a~WW-morphism
$X \leftarrow \mathbf{1}$.
Of course, this is the general case, too, since we may now replace~$X$ by $X \otimes Y$ where~$Y$ is no longer
necessarily the unit object.
In many categories, the ordinary morphisms $X\leftarrow \mathbf{1}$ have a~special meaning.
For instance, in categories of relations, these morphisms correspond to subsets of~$X$; when there is a~further highly
selective structure, such as when the objects are manifolds, the morphisms from $\mathbf{1}$ correspond to submanifolds.
On the other hand, in categories of vector spaces in which $\mathbf{1}$ is the ground f\/ield, the morphisms $X\leftarrow
1$ correspond to the elements of~$X$.
Whenever a~suave morphism $X\leftarrow 1$ in ${\mathbf C}$ has a~name, such as ``thingie'', we will refer to a~diagram
$X\twoheadleftarrow R \leftarrowtail \mathbf{1}$ as a~\textit{hyperthingie}.
Thus, we may speak of ``hypersubsets,'' ``hyperelements'', etc.

\section{Hyperrelations and their trajectories}
In this section, we will work in the category $\mathbf{REL}$, with all morphisms suave and monic pairs congenial.
Everything we do has, though, has a~straightforward extension to other selective categories of relations, including
$\mathbf{SLREL}$, as long as monicity is implied by congeniality.

We will also use the usual rigid monoidal structure in $\mathbf{REL}$, where the monoidal project is the Cartesian
product (denoted by the usual $\times$ rather than $\otimes$), $\mathbf{1}$ is the one-element set $\{\phi\}$, and
${\overline{Y}} = Y$, with the usual unit and counit morphisms.
A~morphism in ${\rm WW}(\mathbf{REL})$ will be called a~\textit{hyperrelation}.
The highly selective structure we will use is the one in which the reductions are the single-valued and surjective
relations; their transposes, the coreductions, are everywhere-def\/ined and injective.
Notice that a~hypergraph $X\times{\overline{Y}} \stackrel{C}{\twoheadleftarrow} R \stackrel{L}{\leftarrowtail}
\mathbf{1}$ is just a~set~$R$ with designated subsets~$C$ and~$L$, along with a~surjection to $X\times {\overline{Y}}$
from~$C$.

If we do not have extra structures to contend with, we can capture much of the essential information in
$X\times{\overline{Y}} \stackrel{C}{\twoheadleftarrow} R \stackrel{L}{\leftarrowtail} \mathbf{1}$ by the induced mapping
$X\times {\overline{Y}} \leftarrow C\cap L$.
Such a~mapping, which is simply a~relation only when it is injective, is sometimes known as a~\textit{span} (of sets).
The spans may be characterized as those diagrams of relations $\xfzgy$ for which~$f$ and~$g$ are everywhere-def\/ined and
single-valued, i.e.~functions.

To generalize from hyperrelations to general composable sequences of morphisms, we introduce the following notion.

\begin{Definition}
A~\textit{trajectory} for a~composable sequence $f=(\ldots, f_{-1},f_0,f_1,\ldots)$ of morphisms in $\mathbf{REL}$ is
a~sequence $(\ldots, x,x,\ldots, x_{-1}, x_0,x_1,\ldots,y,y,\ldots)$ for which $(x_{i-1},x_i)\in f_i$ for all~$i$.
The set of all trajectories for~$f$ will be denoted by ${\cal T}(f)$.
\end{Definition}

When ${\mathbf C}$ is a~selective category of relations in which the suave morphisms have extra structure as objects in
a~category ${\mathbf M}$ admitting f\/ibre products (though the f\/ibre product may no longer represent a~suave morphism),
we may give extra structure to the trajectory spaces as well.
If $f=(f_1,\ldots, f_n)$ is a~path of morphisms $X_{j-1} \stackrel {f_j}{\leftarrow} X_j$ in ${\mathbf C}$, then we may
identify $\tau(f)$ with the f\/ibre product $f_1 \times_{X_1} f_2 \times_{X_2} \dots \times_{X_{n-2}}f_{n-1}
\times_{X_{n-1}} f_n$, so that it is again an object in~$M$ and so shares some of the structure carried by the suave
morphisms themselves.
For example, when~${\mathbf C}$ is a~category of linear relations, the trajectory spaces are themselves linear spaces,
and the source and target maps are linear.

The pairs $(x,y)$ which appear at the ends of trajectories for~$f$ comprise the relation $c'(f)\in \mathbf{REL} (X,Y)$,
where $c'$ is the composition functor def\/ined in Proposition~\ref{prop-pathcomposition}.
We thus have a~map $c'(f)\stackrel{\tau(f)}{\longleftarrow} {\cal T} (f)$ for each path~$f$, and maps $X\leftarrow
\tau(f)$ and $Y\leftarrow \tau(f)$.

\begin{Proposition}
If $f'$ and $f''$ are equivalent paths, then there is is a~bijection ${\cal T}(f')\stackrel{\theta
(f',f'')}{\longleftarrow} {\cal T}(f'')$ for which $\tau(f') \theta(f',f'') = \tau(f'')$.
$($Recall that $c'(f')=c'(f''))$. These bijections may be chosen coherently, i.e., such that
$\theta(f',f'')\theta(f'',f''') = \theta(f',f''')$.
\end{Proposition}
\begin{proof}
For any~$f$, let $f_{m}$ be the minimal representative of $\langle f \rangle$, as described in
Remark~\ref{rmk-minimal}.
We may reduce~$f$ to $f_{m}$ by the following series of steps.
First shift it so that the f\/irst nonidentity entry is~$f_1$.
(If all entries of~$f$ are identities, then~$f$ is already minimal.) Next, eliminate identities from left to right until
there are none left between nonidentity entries.
At each step, there is an corresponding bijection between trajectory sets given either by shifting or by removing
repeating entries.
The composition of these bijections is a~bijection which we def\/ine to be $\theta(f_m,f)$.
Now, for any equivalent $f'$ and $f''$, $f_m'=f_m''$, and we def\/ine $\theta(f',f'')$ to be
$\theta(f_m',f')^{-1}\theta(f_m'',f'')$. Coherence follows immediately.
\end{proof}

More interesting is that, when we pass from $P(\mathbf{REL})$ to ${\rm WW}(\mathbf{REL})$ by using the selective structure
given by monic compositions, the trajectory space of a~morphism $[f]$ and its projection to the relation $c([f])$ are
still well def\/ined.

\begin{Proposition}
If $f'$ and $f''$ are representatives of the same morphism $[f] \in {\rm WW}(\mathbf{REL})$, there is a~bijection ${\cal
T}(f')\stackrel{\theta (f',f'')}{\longleftarrow} {\cal T}(f'')$ for which $\tau(f') \theta(f',f'') = \tau(f'')$.
$($Recall that $c(f')=c(f'').)$ These bijections may be chosen coherently, i.e., such that $\theta(f',f'')\theta(f'',f''')
= \theta(f',f''')$.
\end{Proposition}
\begin{proof}
It is clear that inserting or removing identity morphisms from a~composable~$n$-tuple does not change the set of
trajectories.
This leaves the moves in which two adjacent entries forming a~congenial pair are replaced by their composition.
Suppose that, in a~sequence, these are~$f_j$ and~$f_{j+1}$.
Then any trajectory is a~sequence of the form $(\ldots,x_{j-1},x_j,x_{j+1},\ldots)$, with each~$x_i$ in~$X_i$.
We obtain all the trajectories for the shortened sequence by removing the $x_j$ entries and moving all those to the
right of it one place to the left.
Since the congenial pair $(f_j,f_{j+1})$ is monic, this operation is a~bijective map between spaces of trajectories.
Coherence follows from the fact that a~trajectory has a~def\/inite identity as a~sequence of $x_i$'s, so that a~sequence
of moves beginning and ending with the same sequence of $f_i$'s must act as the identity on the space of trajectories.
\end{proof}

We may therefore refer to the space ${\cal T}([f])$ of trajectories for a~morphism in ${\rm WW}(\mathbf{REL})$, with its
projection $\tau([f])$ to $\mathbf{REL}$.
We can actually reconstruct $[f]$ from its trajectories.

\begin{Proposition}
\label{prop-trajectoriesrepresent}
Any morphism $[\xfy]$ in ${\rm WW}(\mathbf{REL})$ is represented by $X \leftarrow {\cal T}([f])\leftarrow Y$, where the arrows
are the projection to~$X$ and the transpose of the projection to~$Y$.
\end{Proposition}
\begin{proof}
Let $X\times{\overline{Y}} \stackrel{C}{\twoheadleftarrow} R \stackrel{L}{\leftarrowtail} \mathbf{1}$ be a~hypergraph
representing the graph of $[f]$.
The trajectories here may be identif\/ied with the elements of $C \cap L$, which may in turn be identif\/ied with the
trajectories of $[f]$.
We claim that $X\times{\overline{Y}} \stackrel{C}{\twoheadleftarrow} R \stackrel{L}{\leftarrowtail} \mathbf{1}$ is
equivalent to $X\times{\overline{Y}} \stackrel{C'}{\leftarrow} C\cap L \stackrel{L'}{\leftarrowtail} \mathbf{1}$, where
$C'$ and $L'$ are the obvious restrictions of~$C$ and~$L$.
If the inclusion of $C\cap L$ in~$R$ is denoted by~$\iota$, then $C'= C\iota$, and $\iota L' = L$.
The pair $(\iota,L')$ is congenial because~$\iota$ is a~coreduction, and $(C,\iota)$ is congenial because both entries
are single-valued; hence, $[f] = [C,L] = [C,\iota L'] = [C \iota, L'] = [C',L']$.
It follows that $[f]$ itself is represented by $X \leftarrow {\cal T}([f]) \leftarrow Y$.
\end{proof}

We may compose trajectories in the obvious way.
Given $\xfygz$ in ${\rm WW}({\mathbf C})$, we def\/ine ${\cal T}(f){\cal T}(g)$ to be the f\/ibre product ${\cal T}(f)\times_Y
{\cal T}(g)$, the usual composition of spans.
It is easy to see now that ${\cal T}(fg)$ is in natural bijective correspondence with ${\cal T}(f){\cal T}(g)$.
Thus, in the absence of additional structure, the category ${\rm WW}(\mathbf{REL})$ is isomorphic to the category of
isomorphism classes of spans.

\begin{Example}
In the WW-category of smooth canonical relations, any morphism $\mathbf{1} \leftarrow \mathbf{1}$ is given by a~pair of
Lagrangian submanifolds in some manifold~$Q$, so the trace of any $X\leftarrow X$ is essentially a~Lagrangian
intersection.
For instance, the trace of a~morphism given simply by a~Lagrangian submanifold of $X\times {\overline{X}}$ is its
intersection with the diagonal.
For a~morphism whose graph is a~Lagrangian hypersubmanifold
\begin{gather*}
X \stackrel {C}{\twoheadleftarrow}Q\stackrel{L}{\leftarrowtail}\mathbf{1},
\end{gather*}
the inverse image of the diagonal in $X\times {\overline{X}}$ under the projection from~$C$ is a~Lagrangian submanifold
of~$Q$, and its intersection with~$L$ is the trace.

The classif\/ication of endomorphisms in ${\rm WW}(\mathbf{SREL})$ of the one-point manifold is already quite dif\/f\/icult.
If the morphism is given by $L_1$ and $L_2$ in~$Q$ and also by $L_1'$ and $L_2'$ in $Q'$, con\-sideration of trajectories
shows that the intersections $L_1 \cap L_2$ and $L_1'\cap L_2'$ must be isomorphic in various senses (e.g.
they have the same number of elements), but it seems likely that information about the intersection is not enough to
determine the morphism.
The proof of Proposition~\ref{prop-trajectoriesrepresent} fails here because the inclusion of the intersection $L_1 \cap
L_2$ into $Q \times {\overline{Q}}$ is not a~canonical relation.
Exceptionally, perhaps, all the Lagrangian pairs intersecting transversally in a~single point represent the identity
morphism.

We will see in the following section that the classif\/ication of morphisms becomes much simpler if we require our
relations to be linear.
\end{Example}

\section{Linear hypercanonical relations}
In the next few sections of this paper, we will analyze in detail the WW-category built from the highly selective rigid
monoidal category $\mathbf{SLREL}$ of f\/inite dimensional symplectic vector spaces and linear canonical relations.

The monoidal product in $\mathbf{SLREL}$ is the usual set-theoretic Cartesian product, so we will denote it by $\times$
rather than $\otimes$, just as we have been doing for relations in general.
The unit object $\mathbf{1}$ is the zero dimensional vector space whose only element is the empty set.

The dual ${\overline{X}}$ of~$X$ is the same vector space, but with its symplectic structure multiplied by~$-1$.
The morphisms $X\leftarrow Y$ are the Lagrangian subspaces of $X\times {\overline{Y}}$.\footnote{Here, the
identif\/ication of morphisms with their graphs is essentially tautological.} In particular, the morphisms $X\leftarrow
\mathbf{1}$ are just the Lagrangian subspaces of~$X$; similarly, the morphisms $\mathbf{1} \leftarrow X$ are the
Lagrangian subspaces of ${\overline{X}}$, but these are the same as the Lagrangian subspaces of~$X$.
Thus, the unit $\delta_X$ and counit $\epsilon_X$ are both given by the diagonal $\{(x,x) | x \in X\}$, a~Lagrangian
subspace of $X \times {\overline{X}}$.

We will use the selective structure in which all morphisms are suave, but only the monic pairs are congenial, just as in
$\mathbf{REL}$.
Monicity for linear relations $\xfy$ and $\ygz$, def\/ined as injectivity of the projection from $(f\times g)\cap (X
\times \Delta_Y \times {\overline{Z}})$ to $X\times {\overline{Z}}$, is equivalent to injectivity over~$0$, i.e.~the
condition $(f \times g)\cap (\{0_X\} \times \Delta_Y \times \{0_{\overline{Z}}\} = \{0\}$. By elementary symplectic
linear algebra, this is equivalent to $f \times g$ being transversal to $X \times \Delta_Y \times {\overline{Z}}$,
i.e.~transversality of the composition.
For the highly selective structure, we def\/ine the reductions to be those morphisms which are surjective and single
valued and the coreductions those which are injective and everywhere def\/ined, just as in $\mathbf{REL}$.

\begin{Proposition}
The category $\mathbf{SLREL}$ with the structures described above is a~highly selective rigid monoidal category.
\end{Proposition}

\begin{proof}
Most of the required properties follow from those of $\mathbf{REL}$, as demonstrated in
Examp\-les~\ref{ex-selectivesuave},~\ref{ex-relhighlyselective}, and~\ref{ex-relationsmonoidal},
and~\ref{ex-relationsselectiverigidity}.
Morphisms from and to the zero vector space are clearly reductions and coreductions respectively, and application of
Proposition~\ref{prop-factorization} gives the factorization of suave morphisms.
\end{proof}

\begin{Remark}
Many other categories closely related to $\mathbf{REL}$ are also highly selective rigid.
Since monicity implies transversality, the category $\mathbf{SREL}$ of smooth canonical relations is highly selective
rigid.
In fact, the same is true for the category $\mathbf{MREL}$ of all smooth relations.
To see this, it is necessary to prove that the compositions in condition 3 of Def\/inition~\ref{dfn-selectivemonoidal} are
transversal as well as monic.
This is really a~statement about the category $\mathbf {LREL}$ of all linear relations between vector spaces.
It may be checked directly or derived as a~consequence of the fact that the linear duality of vector spaces exchanges
unit and counit morphisms, and monic and transversal compositions.

We also note that the rigid monoidal categories $\mathbf{MREL}$ and $\mathbf{SREL}$ remain highly selective if the
reductions are required to be locally trivial f\/ibrations on their domains, and not just surjective submersions.
In fact, we can take as reductions any rigid monoidal subcategory containing morphisms unit objects and identity
morphisms.
This may be of assistance when one tries to classify hypermorphisms and their hypergraphs, as we will now do for
$\mathbf{SLREL}$.
\end{Remark}

Since the congenial pairs in $\mathbf{SLREL}$ are monic, there is a~good notion of trajectories for the WW-morphisms.
Given $[f]$, if $(x,y)$ is in the shadow of $[f]$, the set of trajectories to~$x$ from~$y$ is an af\/f\/ine space modeled on
the vector space of trajectories to $0$ from~$0$, which is the kernel of a~projection from a~f\/ibre product.
That the dimension of this kernel is an invariant is a~consequence of the invariance of spaces of trajectories, with the
linear structure taken into account.
This justif\/ies the following def\/inition:

\begin{Definition}
Let $\wwxfy$ be a~morphism in ${\rm WW}(\mathbf{SLREL})$.
The \textit{excess} ${\cal E}([f])$ of $[f]$ is the dimension of the af\/f\/ine space of trajectories between any two points
in~$X$ and~$Y$.
\end{Definition}

We omit the easy proof of the following:
\begin{Proposition}
A~pair $(f,g)$ is congenial iff the excess of $[f,g]$ is zero.
More generally, a~WW-morphism is represented by a~single linear canonical relation if and only if its excess is zero.
\end{Proposition}

To classify all of the WW-morphisms, it is useful to begin with the case $Y = \mathbf{1}$.
We def\/ine a~\textit{Lagrangian hypersubspace} of the symplectic vector space~$X$ to be a~WW-morphism $X\leftarrow\mathbf{1}$.
It is represented by diagrams of the form
\begin{gather*}
X \stackrel {C}{\twoheadleftarrow}Q\stackrel{L}{\leftarrowtail}\mathbf{1},
\end{gather*}
where~$C$ is a~reduction whose domain, a~coisotropic subspace, will also be denoted by~$C$, and~$L$ is Lagrangian in~$Q$.
We may therefore denote the Lagrangian hyperspace by $[C,L]$.
The set-theoretic composition $CL$, a~Lagrangian subspace of~$X$, is the shadow of $[C,L]$ and is therefore
well-def\/ined.
The excess of~$[C,L]$ is the dimension of~$C\cap L$ or, equivalently, the codimension of~$C+L$ in~$Q$.

The following fact must be well known, but we could not f\/ind a~proof in the literature.
We will use it to classify Lagrangian hypersubspaces.
(A similar argument is contained in the proof of Theorem~\ref{thm-sabotquotient} below.)

\begin{Lemma}
\label{lemma-isotropiclagrangian}
In a~symplectic vector space~$V$, a~pair consisting of a~Lagrangian subspace~$L$ and an isotropic subspace~$I$ is
characterized up to symplectomorphism by the dimensions of~$V$,~$I$ and $L\cap I$.
\end{Lemma}
\begin{proof}
Choose complements to write $I = (I\cap L) \oplus J$ and $L = (I\cap L)\oplus K$.
Then $J\cap L = \{0\}$, and~$J$ can be extended to a~Lagrangian complement~$M$ of~$L$ in~$V$.
We may identify~$M$ with~$L^*$ and hence identify~$V$ with the direct sum~$L\oplus L^*$ with the usual symplectic
structure.
The decomposition of~$L$ gives, via the orthogonality relation $^{\rm o}$ between subspaces of~$L$ and those~$L^*$, a~dual
decomposition~$L^* = K^{\rm o} \oplus (I\cap L)^{\rm o}$, which is naturally isomorphic to $(I\cap L)^* \oplus K^*$.
Since $J\subset L^*$ is symplectically orthogonal to $(I\cap L)$, it must be contained in the summand~$K^*$.
We may now choose a~basis $(e_i,f_j,g_k)$ of~$L$ with dual basis $(e_i^*,f_j^*,g_k^*)$, forming together a~canonical
basis of~$V$, such that the $e_i$ are a~basis of $I\cap L$ and the $f_j^*$ are a~basis of $J \subseteq K^* \subseteq
L^*$.
The~$e_i$ and~$f_j$ together then form a~basis for~$I$, and this gives a~normal form for the arrangement of spaces
$(V,L,I)$ which depends only on the dimensions of~$V$,~$I$ and $L\cap I$.
\end{proof}

\begin{Proposition}
Two Lagrangian hypersubspaces are equal if they have the same shadow and the same excess.
\end{Proposition}
\begin{proof}
Let~$\Lambda$ be a~Lagrangian subspace of~$X$ and~$k$ a~nonnegative integer.
We will construct a~normal form which is equivalent to any representative of a~Lagrangian hypersubspace with
shadow~$\Lambda$ and excess~$k$.
In $X \times {\mathbb R}^{2k}$, let $C_{k,0}= X \times{\mathbb R}^{k}$ and $\Lambda_{k,0} = \Lambda \times{\mathbb
R}^k$.
Here, ${\mathbb R}^{2k}$ is identif\/ied with the symplectic vector space ${\mathbb R}^k \times {{\mathbb R}^k}^*$.
The Lagrangian hypersubspace $\langle C_{k,0}, \Lambda_{k,0}\rangle$ has shadow~$\Lambda$ and excess~$k$.
This is a~minimal representative in its equivalence class; we get larger representatives by forming its monoidal product
with the trivial Lagrangian hypersubspace of the point ${\mathbb R}^0$, as represented as the transversal pair
$({\mathbb R}^r,{{\mathbb R}^r}^*)$ with intermediate space ${\mathbb R}^{2r}$.
Denote this product by $\langle C_{k,r}, \Lambda_{k,r}\rangle$.

Now let $\langle C,L \rangle$ be any Lagrangian hypersubspace of~$X$ with shadow~$\Lambda$, excess~$k$, and intermediate
space~$Q$ of dimension $2N$.
If~$X$ has dimension $2n$, then the dimension of~$C$ must be $N+n$.
We observe f\/irst that the diagram
\begin{gather*}
X \stackrel {C}{\twoheadleftarrow}Q\stackrel{L}{\leftarrowtail} \mathbf 1
\end{gather*}
is symplectically isomorphic to
\begin{gather*}
X \stackrel {C_{k,r}}{\twoheadleftarrow}X\times{\mathbb R}^{2k}\times{\mathbb
R}^{2r}\stackrel{\Lambda_{k,r}}{\leftarrowtail} \mathbf 1;
\end{gather*}
i.e.~there is a~symplectomorphism of~$Q$ with $X \times {\mathbb R}^{2k} \times {\mathbb R}^{2r}$ with $r=N-n-k$,
taking~$L$ to $\Lambda \times {\mathbb R}^k \times{{\mathbb R}^r}^*$ and~$C$ to $X\times {\mathbb R}^k\times {{\mathbb
R}^r}$.
The isomorphism now follows from Lemma~\ref{lemma-isotropiclagrangian} above.
That $\langle C,L\rangle$ and $\langle C_{k,r}, \Lambda_{k,r}\rangle$ are equal as WW-morphisms now follows from
Proposition~\ref{prop-hyperisomo}.
\end{proof}

Thus, there is a~bijective correspondence between Lagrangian hypersubspaces of~$X$ and pairs $(L,k)$, where~$L$ is an
ordinary Lagrangian subspace and~$k$ is a~nonnegative integer.

We can now understand general WW-morphisms via their graphs.
We will call any morphism $\wwxfy$ in ${\rm WW}(\mathbf{SLREL})$ a~\textit{(linear) hypercanonical relation} to~$X$ from~$Y$.
Its graph is a~Lagrangian hypersubspace of $X\times {\overline{Y}}$.
These have the following two useful properties.

\begin{Proposition}\quad
\begin{enumerate}\itemsep=0pt
\item[$(1)$] The excess of any hypercanonical relation is equal to that of its graph.
\item[$(2)$] If $[f]$ and $[g]$ are composable hypercanonical relations,
\begin{gather*}
{\cal E}([f][g]) = {\cal E}([f]) + {\cal E}([g]) + {\cal E}([c([f]),(c[g])]).
\end{gather*}
\end{enumerate}
\end{Proposition}
\begin{proof}
Given a~morphism $[f]$ represented by $X \stackrel{a}{\twoheadleftarrow} Q \stackrel{b}{\leftarrowtail} Y$, its excess
is the dimension of the trajectory space
\begin{gather*}
(a \times b)\cap (\{0_X\}\times \Delta_Q\times \{0_Y\}) \subseteq X\times Q\times {\overline{Q}}\times {\overline{Y}}.
\end{gather*}

On the other hand, the graph is represented~by
\begin{gather*}
X\times {\overline{Y}} \stackrel{a \times 1_{\overline{Y}}}{\twoheadleftarrow} Q \times {\overline{Y}}
\stackrel{\gamma_b}{\leftarrowtail} \mathbf{1}.
\end{gather*}
The excess of the latter is the dimension of
\begin{gather*}
(a\times 1_{\overline{Y}} \times \gamma_b) \cap({0_{X\times {\overline{Y}}}} \times \Delta_{{\overline{Q}} \times Y}
\times\{0_\mathbf{1}\})\subseteq X \times {\overline{Y}}\times {\overline{Q}}\times Y \times Q\times{\overline{Y}}
\times \mathbf{1}.
\end{gather*}
This intersection consists of the sextuples $(x,y,q',y',q'',y'')$ such that $x=0_X$, $y=0_Y$, $y=y'=y''$, $q'=q''$,
$(x,q')\in a$, and $(q'',y'')\in b$. These may be identif\/ied with the trajectories to $0_X$ from $0_Y$ in $[a,b]=[f]$.

Note that the graph of $[f]$ is also represented by the ``graph product''
\begin{gather*}
X\times {\overline{Y}} \twoheadleftarrow X \times {\overline{Q}} \times Q \times {\overline{Y}} \stackrel{\gamma_a
\times \gamma_b}{\leftarrowtail} \mathbf{1}\times\mathbf{1}=\mathbf{1},
\end{gather*}
whose trajectory space is essentially the same as that of $[ab]$.

For (2), we begin with the fact that ${\cal T}([f,g])$ is the f\/ibre product over~$Y$ of ${\cal T}([f])$ and ${\cal
T}([g])$.
The projection $\tau([f,g])$ may therefore be factored as
\begin{gather*}
c([fg])\stackrel{}{\twoheadleftarrow} c([f])\times_Y c([g]) \stackrel{}{\leftarrowtail} {\cal T}([f,g]).
\end{gather*}
The kernel of map to $c([fg])$ has dimension ${\cal E}([c([f]),(c[g])])$, while the kernel of the map from ${\cal
T}([f,g])$ has dimension ${\cal E}([f]) + {\cal E}([g])$.
\end{proof}

We may therefore identify the WW-morphisms (i.e.~hypercanonical relations) $X\leftarrow Y$ with the pairs $(f,k)$,
where~$f$ is a~Lagrangian subspace of $X \times {\overline{Y}}$ and~$k$ is a~nonnegative integer.
We will call these pairs \textit{indexed canonical relations}.
When~$Y$ is the unit object, we will call the pairs \textit{indexed Lagrangian subspaces} of~$X$.

The following theorem expresses the structure of ${\rm WW}(\mathbf{SLREL})$ via the identif\/ication with indexed canonical
relations.
The proof consists of elementary dimension calculations.

\begin{Theorem}
The indexed $($linear$)$ canonical relations form a~category by identification with the linear hyperrelations.
The composition law is
\begin{gather*}
(f',k')(f'',k'') = (f'f'',k'+k''+{\cal E}([f',f'']),
\end{gather*}
and the monoidal product is
\begin{gather*}
(f',k')\times(f'',k'') = (f'\times f'',k'+k'').
\end{gather*}
The monoid of endomorphisms of the unit object is naturally identified with the nonnegative integers, and its action on
the category is the shifting operation $k'\cdot (f,k)= (f,k'+k)$. The trace of an endomorphism $(f,k)$ of~$X$ is $ \dim
(\gamma_f \cap \Delta_X)+k$.
\end{Theorem}

Note that $\gamma_f \cap \Delta_X$ is (the diagonal of) the ``f\/ixed point space'' of~$f$.

\section{Indexed Lagrangian Grassmannians and the Sabot topology}
We def\/ine the \textit{indexed Lagrangian Grassmannian} of a~symplectic vector space~$X$ to be the product ${\cal
L}_{\bullet} (X) $ of the usual Lagrangian Grassmannian ${\cal L}(X)$ with the nonnegative integers.
It is the union of ``levels'' ${\cal L}_k (X) = {\cal L} (X) \times \{k\}$, with~$k$ called the \textit{indExample} of
$(L,k)$.
This space carries a~useful topology weaker than the usual one, rendering the operations of composition and reduction
continuous.

To def\/ine our topology, we we f\/irst introduce a~(discrete) metric\footnote{This metric may already by found in the
paper~\cite{ch:geometry} by Chow, who proves that isometries must come from projective transformations.}~$d$ on ${\cal
L}(X)$, with $d(L,L')$ being the codimension ${\rm codim\,}(L \cap L',L)$ of $L\cap L'$ in~$L$ (or in $L'$).
We then def\/ine the {\em sublevel} operator~$S$ from points to subsets in ${\cal L}_{\bullet} (X) $~by
\begin{gather*}
S(L,k) = \{(L',j) | d(L,L')\leq k-j\}.
\end{gather*}
The triangle inequality for~$d$ implies that the sublevel operator def\/ines a~partial order on the extended Lagrangian
Grassmannian, def\/ined by $(L,k)\leq (L',k')$ when $(L,k) \in S(L',k')$. In other words,
\begin{gather*}
(L,k)\leq (L',k')
\qquad
\text{if\/f}
\qquad
d(L,L') \leq k'-k.
\end{gather*}
The partial ordering property means in particular that $S\circ S = S$.

\begin{Lemma}
The subsets $S({\cal U})$, where ${\cal U}$ ranges over the open subsets of the levels in ${\cal L}_\bullet (X)$, form
the basis of a~topology.
\end{Lemma}

\begin{proof}
For $j=1,2$, let ${\cal U}_j$ be an open subset of ${\cal L}_{n_j}(X)$, and let $(L,k)$ belong to $S({\cal U}_1) \cap
S({\cal U}_2)$.
We must f\/ind a~neighborhood ${\cal U}$ of $(L,k)$ such that $S({\cal U})$ is contained in that intersection.

Let ${\cal U}$ be the intersection of $S({\cal U}_1) \cap S({\cal U}_2)$ with ${\cal L}_k (X)$.
Then $(L,k)\in {\cal U}$, and
\begin{gather*}
S {\cal U} \subseteq S(S({\cal U}_1))\cap S(S({\cal U}_2))= S({\cal U}_1) \cap S({\cal U}_2).
\tag*{\qed}
\end{gather*}
\renewcommand{\qed}{}
\end{proof}

We call the topology generated by the sublevels of open sets in the levels the \textit{Sabot topology}.
We choose this name because, as we are about to see, a~result in Sabot~\cite{sa:electrical} implies that this topology
renders an indexed symplectic reduction operation continuous.
Note that the Sabot topology satisf\/ies the separation axiom~$T_0$, but not~$T_1$.
In particular, it is not Hausdorf\/f.
We also observe that each shift operation taking~$(L,k)$ to~$(L,k+k')$ is a~homeomorphism from~${\cal L}_\bullet (X)$
onto its image.

Although the primary operation in our category is composition, it is useful to analyze it via the simpler operation of
reduction.

Let~$C$ be a~coisotropic subspace of~$X$, $X^C = C/C^\perp$the reduction of~$X$ by~$C$.
Associated with this data is a~canonical relation $X^C \stackrel{\rho^C}{\twoheadleftarrow} X$, the graph of the
projection $C\leftarrow C^\perp$.
Identifying Lagrangian subspaces $L\subset X$ with morphisms $X\leftarrow \mathbf{1}$, we may compose them with
$\rho^C$, obtaining a~map ${\cal L} (X^C)\leftarrow {\cal L} X$ which we will also denote by $\rho^C$.
In terms of subspaces, it takes each~$L$ to $ (L\cap C)/(L\cap C^\perp)$, which will denote by $L^C$.
Although this operation on Lagrangian subspaces is a~very natural one, it is discontinuous at all~$L$ for which~$L$ is
not transversal to~$C$.

To quantify (and correct) this discontinuity, we begin by observing that the codimension of $L+C$ in~$X$ is equal to the
dimension of $L\cap C^\perp$, which is precisely the excess ${\cal E}(\rho^C,L)$ of the pair consisting of the reduction
and inclusion morphisms.
We will denote it by the abbreviated notation ${\cal E}^C(L)$.
In terms of this excess, we def\/ine the \textit{indexed} reduction operation~$\rho_\bullet^C$ to~${\cal L}_\bullet (X^C)$
from~${\cal L}_\bullet (X)$ by $\rho_\bullet^C (L,k) = (L^C,k+{\cal E}^C(L))$.
It corresponds precisely to the operation of composition in the category ${\rm WW}(\mathbf{SLREL})$ with the canonical
relation~$\rho^C$, which takes WW-morphisms $X\leftarrow \mathbf{1}$ to morphisms~$X^C \leftarrow \mathbf{1}$.

\begin{Remark}
As is well known, we can also express composition in terms of reduction.
Namely, given~$X$,~$Y$, and~$Z$, there is a~natural identif\/ication of $X \times {\overline{Z}}$ with the reduction of $X
\times {\overline{Y}} \times Y \times {\overline{Z}}$ by $X\times \Delta_Y \times {\overline{Z}}$.
The composition $\xfygz$ is then the reduction of the monoidal (i.e.~Cartesian) product $f \times g$.
In this way, we can recover the formula for indexed composition from that for indexed reduction (along with additivity
with respect to monoidal product).

We also note that reduction commutes with the shift operation.
\end{Remark}

Now we will prove the continuity of symplectic reduction, using a~ref\/inement of the following theorem of Sabot.

\begin{Theorem}
[\cite{sa:electrical}] Let~$C$ be a~coisotropic subspace of~$X$.
The closure $\overline \rho^C \subset {\cal L} (X) \times {\cal L} (X^C)$ of the graph of $\rho^C$ is the multivalued
function whose graph consists of those pairs $(L,\Lambda)$ for which $d(L^C,\Lambda)\leq {\cal E}^C(L)$.
\end{Theorem}

Sabot also proves that $\overline \rho^C$ is a~rational variety, i.e.~a ``rational map'' (in fact, this is what is
really important for his purposes), but we will not use this result here.

Our ref\/inement goes as follows.

\begin{Theorem}
\label{thm-refinedsabot}
Let~$C$ be a~coisotropic subspace of~$X$, $L_i$ a~sequence of Lagrangian subspaces in~$X$ which converges to
a~limit~$L$, with $L^C_i$ convergent to $L'\in {\cal L}(X^C)$.
Then
\begin{gather*}
d(L^C, L') \leq {\cal E}^C (L) - \limsup_{i\to \infty} {\cal E}^C (L_i).
\end{gather*}
\end{Theorem}

\begin{proof}
By passing to a~subsequence, we may assume that ${\cal E}^C(L_i)$ is constant and equal to the $\limsup$ in the
statement, which we will denote simply by ${\cal E}$.
We also assume that $L_i \cap C$ converges to a~subspace~$R$ and that $L_i \cap C^\perp$ converges to a~subspace~$S$.
If~$X$ has dimension~$2n$ and~$C$ has dimension $n+m$ (so that $X^C$ has dimension $2m$), then~$R$ has dimension $m +
{\cal E}$ and~$S$ has dimension~${\cal E}$.
The convergence of $L_i$ to~$L$ implies that $R \subseteq L \cap C$ and $S \subseteq L \cap C^\perp$. The inclusions
above are not necessarily equalities, but they induce a~map from $R/S$ to $X^C$.
The image $R^C = R/(R\cap C^\perp)$ is isotropic and, although $R/S$ has half the dimension of~$X$, this map is not
necessarily an embedding.
Since~$R$ is contained in~$L$, $R^C$ is contained in $L^C$, and the convergence implies that it is contained in $L'$ as
well, hence in $L^C \cap L'$.
The dimension of the image is determined by that of $R\cap C^\perp$.
Since the latter is contained in $L\cap C^\perp$, its dimension is at most ${\cal E}^C(L)$, so
\begin{gather*}
d(L^C,L') \leq m - (\dim (R) - \dim (R \cap C^\perp)) \leq m - (m + {\cal E} - {\cal E}^C(L) = {\cal E}^C(L) - {\cal E}.
\tag*{\qed}
\end{gather*}
\renewcommand{\qed}{}
\end{proof}

\begin{Remark}
One can also show that all values of the distance between $0$ and ${\cal E}^C (L) - \limsup\limits_{i\to \infty} {\cal E}^C
(L_i)$ are realized, but we will not use this fact.
\end{Remark}

We will also use the following lemma to the ef\/fect that reduction is distance non-increasing.

\begin{Lemma}
\label{lemma-reduction}
Let~$C$ be a~coisotropic subspace of~$X$, $L_1$ and $L_2$ Lagrangian subspaces.
Then $d(L^C_1,L^C_2) \leq d(L_1,L_2)$.
\end{Lemma}

\begin{proof}
We f\/irst note that, if~$W$,~$Y$, and~$Z$ are three subspaces of some vector space, with $W\subseteq Y$, then ${\rm
codim\,}(W\cap Z,Y\cap Z) \leq {\rm codim\,} (W,Y)$. This follows from the exact sequence
\begin{gather*}
0 \to (Y\cap Z)/(W\cap Z) \to Y/W \to (Y + Z)/(W+Z) \to 0
\end{gather*}
(which also evaluates the dif\/ference between the two codimensions).

Applying the inequality above with $W = L_1 \cap L_2$, $Y = L_1$, and $Z = C$, we f\/ind
\begin{gather*}
{\rm codim\,}((L_1 \cap C)\cap (L_2 \cap C), L_1 \cap C) = {\rm codim\,} (L_1 \cap L_2 \cap C,L_1 \cap C)
\\
\phantom{{\rm codim\,}((L_1 \cap C)\cap (L_2 \cap C), L_1 \cap C)}
\leq {\rm codim\,}(L_1\cap L_2,L_1)
= d(L_1,L_2).
\end{gather*}
Now $L^C_1 \cap L^C_2 $ is the image of $(L_1 \cap C)\cap (L_2 \cap C)$ under the projection along $C^\perp$, and the
intersection of $(L_1 \cap C)\cap (L_2 \cap C)$ with $C^\perp$ is contained in the intersection of $L_1 \cap C$ with
$C^\perp)$, so the codimension of one in the other can only decrease under the projection.
Hence
\begin{gather*}
d\big(L_1^C,L_2^C\big) = {\rm codim\,} (L_1^C\cap L_2^C, L_1^C) \leq d(L_1,L_2).
\tag*{\qed}
\end{gather*}
\renewcommand{\qed}{}
\end{proof}

In fact, there is a~stronger version of the result.

\begin{Proposition}
\label{prop-reduction}
Let~$C$ be a~coisotropic subspace of~$X$, $L_1$ and $L_2$ Lagrangian subspaces.
Then the inequality
\begin{gather}
\label{eq-inequality}
d\big(L_1^C,L_2^C\big) \leq d(L_1,L_2) - \big|{\cal E}^C(L_1) - {\cal E}^C(L_2)\big|.
\end{gather}
is equivalent to the statement that the indexed reduction operation $\rho_\bullet^C$ between indexed Grassmannians is
order-preserving.
These two equivalent statements are true.
\end{Proposition}
\begin{proof}
We begin by proving the equivalence of~\eqref{eq-inequality} with the order-preserving property of indexed reduction.
Suppose that~\eqref{eq-inequality} is true, and suppose that $(L_1,k_1) \leq (L_2,k_2)$. Then $d(L_1,L_2) \leq k_2 -
k_1$, and so the inequality in the statement gives:
\begin{gather*}
d\big(L_1^C,L_2^C\big) \leq k_2 - k_1 + {\cal E}^C(L_2) - {\cal E}^C(L_1),
\end{gather*}
which states precisely that $\rho_\bullet^C(L_1,k_1) \leq \rho_\bullet^C(L_2,k_2)$.
Conversely, given any $L_1$ and $L_2$, we may assume that ${\cal E}^C(L_1)\leq {\cal E}^C(L_2)$.
Assuming order preservation and applying it to the relation $(L_1,d(L_1,L_2)) \geq (L_2,0)$, we
obtain~\eqref{eq-inequality}.

We now prove that these properties are true.
It suf\/f\/ices to prove~\eqref{eq-inequality} for reduction by a~codimension-one subspace, since any reduction is
a~composition of these, and the composition of two order-preserving maps is order-preserving.
When~$C$ has codimension~$1$, then the excess of a~Lagrangian subspace is either $0$ or~$1$.
The equality then follows immediately from Lemma~\ref{lemma-reduction} unless ${\cal E}^C(L_1)=1$ and ${\cal E}^C(L_2) =
0$.
In that case, we have $C^\perp \subseteq L_1\subseteq C$.
Since $L_2$ is independent of $C^\perp$ and transversal to~$C$, it follows that the projection $C \to X^C$ is an
isomorphism when restricted to $C\cap L_2$.
In particular, it is an isomorphism on $L_1 \cap L_2$, so the dimension of intersection remains the same, and hence the
codimension increases by~$1$, consistent with~\eqref{eq-inequality}, which is now proven.
\end{proof}

Here is our continuity result.

\begin{Theorem}
\label{thm-reduction}
Let~$C$ be a~coisotropic subspace of~$X$.
The indexed reduction operation $\rho_\bullet^C$ is continuous for the Sabot topologies on $ {\cal L}_\bullet(X)$ and
${\cal L}_\bullet (X^C)$.
\end{Theorem}

\begin{proof}
If $\rho_\bullet^C$ were not continuous, there would be an $(L,k)$ in ${\cal L}_\bullet(X)$ and an open neighbor\-hood~${\cal V}$ of $L^C$ such that no neighborhood of $(L,k)$ could map entirely into the neighborhood $S({\cal V},k+{\cal
E}^C(L))$ of $\rho_\bullet^C(L,k) = (\rho_C(L),k+{\cal E}_C(L))$.
Choosing a~sequence of neighborhoods of~$(L,k)$ shrinking to the intersection $S(L,k)$ of all such neighborhoods (e.g.~a
sequence $S({\cal U}_i,k)$, where the~${\cal U}_i$ are neighborhoods shrinking down to~$L$), we could f\/ind a~sequence
$(L_i,k_i)$ converging to $(L,k)$ such that $\rho_\bullet^C (L_i,k_i)$ remains outside $S({\cal V},k+{\cal E}^C(L))$.
Passing to a~subsequence if necessary, we could assume: (1) $k_i$ is a~constant, which we will call $k_0$; (2) ${\cal
E}^C(L_i)$ is a~constant, which we will call ${\cal E}$; (3) $L_i$ converges to some $L_0$, which must be in $S(L,k)$,
i.e.~$d(L,L_0) \leq k-k_0$; (4) $\rho^C(L_i)$ converges to some $L'$ such that $(L',k_0 + {\cal E})$ is not contained in
$S({\cal V},k+{\cal E}^C(L))$, in particular, not in $S(L^C,k+{\cal E}^C(L))$.
To contradict this, it suf\/f\/ices to prove that
\begin{gather*}
d(L^C,L') \leq k + {\cal E}^C(L) - (k_0 + {\cal E}).
\end{gather*}

Now $d(L^C_0,L') \leq {\cal E}^C(L_0) - {\cal E}$ by Theorem~\ref{thm-refinedsabot}, $d(L^C,L^C_0) \leq d(L,L_0)-|{\cal
E}^C (L) - {\cal E}^C (L_0)|$ by Proposition~\ref{prop-reduction}, and $d(L,L_0) \leq k-k_0$ by (3) above, so, by the
triangle inequality, $d(L^C,L') \leq k-k_0 - |{\cal E}^C (L) - {\cal E}^C (L_0)|+ {\cal E}^C(L_0) - {\cal E}$. Since $
{\cal E}^C(L_0) - |{\cal E}^C (L) - {\cal E}^C (L_0)| \leq {\cal E}^C(L)$, the proof is complete.
\end{proof}

We turn now to general indexed canonical relations, identif\/ied with indexed Lagrangian subspaces of spaces of the form
$X \times {\overline{Y}}$.
Given $X \stackrel{(f',k')}{\longleftarrow} Y$ and $Y \stackrel{(f'',k'')}{\longleftarrow} Z$, their composition
$(f',k')(f'',k'')$ has been seen to be $(f'f'',k'+k''+{\cal E}^C(f' \times f''))$, where $f'f''\in X \times
{\overline{Z}}$ is the set-theoretic composition, and~$C$ is the coisotropic subspace $X \times \Delta_Y \times
{\overline{Z}}$ in $X \times {\overline{Y}} \times Y \times {\overline{Z}}$. In other words, up to the identif\/ication of
the reduced space $(X \times {\overline{Y}} \times Y \times {\overline{Z}})^C$ with $X \times {\overline{Z}}$, the
composition is the indexed reduction by~$C$ of the monoidal product $f'\times f''$.

Since the composition operation is itself a~composition of the continuous operation of reduction with the monoidal
product, and the monoidal product is easily seen to be continuous, we have:
\begin{Theorem}
\label{thm-composition}
Let~$X$,~$Y$, and~$Z$ be symplectic vector spaces.
Then the indexed composition operation from ${\cal L}_\bullet(X\times {\overline{Y}}) \times {\cal L}_\bullet (Y \times
{\overline{Z}})$ to ${\cal L}_\bullet (X \times {\overline{Z}})$ is continuous in the Sabot topologies.
\end{Theorem}
\begin{Remark}
If we had an independent proof of Theorem~\ref{thm-composition}, Theorem~\ref{thm-reduction} would be an easy corollary,
given the description of reduction as composition with a~particular canonical relation to $X^C$ from~$X$.
\end{Remark}

\section{The Sabot topology as a~quotient topology}
Let~$C$ be a~coisotropic subpace of~$X$, with reduced space $X^C$.
We will consider here the restricted reduction operation, still to be denoted by $\rho_\bullet^C$, from the {\em
ordinary} Lagrangian Grassman\-nian~${\cal L}(X)$ to the truncated indexed Lagrangian Grassmannian ${\cal L}_{\leq
r}(X^C)$ consisting of the levels up to the dimension $r=n-m$ of $C^\perp$.

\begin{Theorem}
\label{thm-sabotquotient}
The restricted Sabot topology on ${\cal L}_{\leq r}X_C$ is the quotient topology for the restricted reduction operation
$\rho_\bullet^C:{\cal L}(X)\rightarrow {\cal L}_{\leq r}(X^C)$ and the usual topology on ${\cal L} (X)$.
\end{Theorem}

\begin{proof}
We already know from Theorem~\ref{thm-reduction} that this reduction map is continuous, so it remains to show that it is
surjective and open.

Given coisotropic~$C$, we may split~$X$ symplectically as $J \oplus (C^\perp \oplus D)$, where~$J$ is any complement of
$C^\perp$ in~$C$ (hence symplectic and projecting isomorphically to $X^C$), and~$D$ is an isotropic complement to~$C$ in~$X$.
$J$ may be identif\/ied with $X_C$ via the projection.
As subspaces of the symplectic subspace $C^\perp \oplus D$, $C^\perp$ and~$D$ give a~Lagrangian splitting, while
$C=C^\perp \oplus J$.

Given now an indexed Lagrangian subspace $(\Lambda,k)$ in $X^C$ with $k \leq \dim C^\perp$, we choose
any~$k$-dimensional subspace $Y \subset C^\perp$ and let $W\subset J$ be the inverse image of~$\Lambda$ under the
projection~$C\rightarrow X^C$.
It may be identif\/ied with~$\Lambda$ itself under the above-mentioned identif\/ication between~$J$ and~$X^C$.
Let $Z\subset D$ be the annihilator of~$Y$ under the symplectic pairing between~$C^\perp$ and~$D$.
Then $L=W\oplus Y \oplus Z$ is Lagrangian, and its indexed reduction~$L^C$ is~$(\Lambda,k)$.
This shows that~$\rho_\bullet^C$ is surjective.

For openness, we must show that, for any Lagrangian subspace $L \subset X$ and any neighbor\-hood~${\cal U}$ of~$L$ in
${\cal L}(X$), the image $\rho_\bullet^C({\cal U})$ contains a~neighborhood of $\rho_\bullet^C(L) = (L^C,k)$, where $k
={\cal E}^C(L)$.

Given~$L$, we f\/irst choose a~decomposition of~$X$ as above, but adapted to~$L$.
First, we decompose $L\cap C$ as a~direct sum $W\oplus Y$, where $Y = L\cap C^\perp$ and $W\subset C$ projects
bijectively to the reduced Lagrangian subspace $L^C \subset X^C$.
We may then choose the complement~$J$ of $C^\perp$ in~$C$ to contain~$W$ (as a~Lagrangian subspace).
$L$ itself may be written as $W \oplus Y \oplus Z$, where~$Z$ is an isotropic subspace of~$X$ projecting injectively to $X/C$.
Now extend~$Z$ to an isotropic complement~$D$ of~$C$ in~$X$.
Then~$Z$ is the annihilator of~$Y$ with respect to the symplectic pairing between $C^\perp$ and~$D$, and the splitting
$W \oplus (Y\oplus Z)$ of~$L$ is compatible with the symplectic decomposition $J \oplus (C^\perp \oplus D)$ of~$X$.

Since~$J$ is identif\/ied with $X^C$, we may f\/ind in $\rho_\bullet^C({\cal U})$ any Lagrangian subspace in a~neighborhood
of $\rho_\bullet^C(L)$ on the same level~$k$ as $\rho_\bullet^C(L)$ simply by moving~$W$ within~$J$ and leaving the
summand $Y\oplus Z\subset C^\perp\oplus D$ f\/ixed.

Next, we will show how to move~$L$ so that $\rho_\bullet^C(L)$ moves to any chosen $(\Lambda,k-q)$ in the sublevel set
$S(\rho_\bullet^C(L))$.
This will complete the proof, since the argument in the previous paragraph may be applied with~$L$ replaced by any
element of ${\cal U}$.

Suppose then, that $d(L^C,\Lambda) = r \leq k -q$ for some~$q$.
To move~$L$ so that its reduction becomes~$\Lambda$ and its level drops by~$q$, we must take an~$r$-dimensional subspace
out of~$W$ (identif\/ied with $L^C$) and replace it by a~dif\/ferent one, while taking a~$q$-dimensional subspace out of $Y
= L \cap C^\perp$.

At this point, it is convenient to choose a~basis of~$X$ suitable for our needs.
Let $2m$ be the dimension of $X^C$ (and hence of~$J$), and let $W'$ be the inverse image in~$J$ of~$\Lambda$.
$W$ is Lagrangian in~$J$.
We start with a~basis $(e_1,\ldots,e_{r},f_1,\ldots,f_{m-r})$ of~$W$ for which the last $m-r$ entries are a~basis of the
subspace $W \cap W'$, which does not have to move.
Now, in the symplectic vector space $(W + W')/(W\cap W')$, the images of~$W$ and $W'$ are complementary Lagrangian
subspaces, so we may take the basis of $W'/(W\cap W')$ which is dual to the images of $e_1,\ldots,e_{r}$ and lift it to
vectors $e_1^*,\ldots,e_r^*$ in $W'$ such that $(e_1^*,\ldots,e_{r}^*,f_1,\ldots,f_{m-r})$ is a~basis of $W'$.
We may then extend this to a~canonical basis
\begin{gather*}
(e_1,\ldots,e_{r},f_1,\ldots,f_{m-r},e^*_1,\ldots,e^*_{r},f^*_1,\ldots,f^*_{m-r})
\end{gather*}
of~$J$.
Next, we choose a~basis
\begin{gather*}
(g_1,\ldots,g_r,h_1,\ldots,h_{q-r},i_1,\ldots,i_{k-q},j_1,\ldots,j_{n-m-k})
\end{gather*}
of $C^\perp$ such that the~$g$'s,~$h$'s, and~$i$'s form a~basis of $Y= L\cap C^\perp$.
Via the pairing with $C^\perp$, we get a~dual basis
\begin{gather*}
(g^*_1,\ldots,g^*_r,h^*_1,\ldots,h^*_{q-r},i^*_1,\ldots,i^*_{k-q},j^*_1,\ldots,j^*_{n-m-k})
\end{gather*}
of~$D$ such that the $j^*$'s form a~basis of the annihilator~$Z$ of~$Y$.
Putting these all together (with boldface a~reminder that indices have been omitted for simplicity of notation), we get
a~canonical basis $({\bf e},{\bf f},{\bf g},{\bf h},{\bf i},{\bf j},{\bf e}^*,{\bf f}^*,{\bf g}^*,{\bf h}^*,{\bf i}^*,{\bf j}^*)$ of~$X$ for which $({\bf e},{\bf f},{\bf g},{\bf h},{\bf i},{\bf j})$ is a~basis of~$L$ and
$({\bf e},{\bf f},{\bf e}^*,{\bf f}^*,{\bf g},{\bf h},{\bf i},{\bf j})$ is a~basis of~$C$, with $({\bf g},{\bf h},{\bf i},{\bf j})$ spanning $C^\perp$.

We are ready to deform~$L$ to a~family $L_t$, for a~real parameter near~$0$, in such a~way that $L_0=L$, while
$\rho_\bullet^C(L_t) = (\Lambda,k-q)$ for all nonzero~$t$.
To do this, we simply deform the basis $({\bf e},{\bf f},{\bf g}, {\bf h},{\bf i},{\bf j}^*)$ of~$L$ to
\begin{gather*}
({\bf e}+t{\bf g}^*,{\bf f},{\bf g}+t{\bf e}^*,{\bf h}+t{\bf h}^*,{\bf i},{\bf j}^*),  
\end{gather*}
which are still linearly independent and thus form a~basis for a~deformed space $L_t$.
One may check by taking all possible symplectic pairings that $L_t$ is still Lagrangian.
When $t=0$, the basis is undeformed.

Suppose now that $t\neq 0$.
Then $L_t \cap C$ has basis $({\bf f},{\bf g}+t {\bf e}^*,i)$, and the basis of $L_t\cap C^\perp$ is simply $\bf i$.  
Modulo $C^\perp$, $L_t $ is the same as $W'$, so $L_t^C = \Lambda$.
Since the basis of $L_t \cap C$ has shrunk from $({\bf g},{\bf h},{\bf i})$ to $\bf i$, the excess has dropped from~$k$ to $k-q$, as
required.
\end{proof}

\section{Negative indices}
Although nonnegative indices are enough to parametrize the linear $\rm WW$-category, there are good reasons to allow
arbitrary integer indices.
This section is an informal discussion of some questions related to the notion of negative indices.

If~$C$ is coisotropic in~$X$, there is an important endomorphism $R^C$ of~$X$ consisting of those pairs $(x,y)$ in
$C\times C$ for which $x-y$ belongs to $C^\perp$.
This relation is the composition\footnote{Note that the excess of this composition is zero.} $(\rho^C)^t\rho^C$; after
quantization, it should become a~projection operator from the vector space quantizing~$X$ to the subspace quantizing
$X^C$.
Since projection operators are idempotent, $R^C$ should be idempotent as well.
This is indeed true for the usual composition, but not for indexed composition.
The square of the indexed canonical relation $(R^C,k)$ is $(R^C,2k + \dim C^\perp)$.
For this to be equal to~$(R^C,k)$ when~$C$ is not all of~$X$ (in which case $R^C$ is the identity and~$k$ may be taken
to be~$0$), we must set~$k$ equal to the negative integer $-\dim C^\perp$.
Since the trace of~$R^C$ is $\dim C$, the trace of $(R^C,-\dim C^\perp)$ is $\dim C - \dim C^\perp$, which is the
dimension (and hence the trace of the identity morphism of) the reduced space $C/C^\perp$.

This example also leads to the introduction of an ``indexed transpose''.
We observe that the set-theoretic composition $\rho^C (\rho^C)^t$ is the identity $1_{C/C^\perp}$, but the excess of the
composition is equal to $\dim C^\perp$.
This suggests that we attach the index $-\dim C^\perp$ to one of the two factors; we propose attaching it to the
multiple-valued relation $(\rho^C)^t$.

Another reason that it seems useful to introduce negative indices is to invert relations $(f,k)$ where $\xfy$ is
invertible; the inverse would be $(f,-k)$.
In particular, we may do this when $f$ is any indexed endomorphism of $\mathbf{1}$.
In fact, these endomorphisms constitute a~multiplicative system, and localizing the $\rm WW$-category there is essentially
the same as allowing negative indices\footnote{In fact, it suf\/f\/ices to localize at the multiplicative system consisting
of the endomorphisms of the unit object.}.

This extension raises many issues.
\begin{itemize}\itemsep=0pt
\item Is there a~version of the ``negative indices'' construction for more general WW-categories? Perhaps one should
localize at the multiplicative system consisting of the endomorphisms of the unit object whose shadows are invertible.
Or maybe, when it does not follow from this, to invert all morphisms whose shadows are suave isomorphisms.

\item Could one introduce some kind of model category structure on ${\rm WW}({\mathbf C})$ when~${\mathbf C}$ is highly selective, with
reductions and coreductions as the f\/ibrations and cof\/ibrations? And then pass to a~homotopy category? What would the
weak equivalences be? Maybe those for which the shadow is a~suave isomorphism.

\item What happens on the quantum side? For instance, consider the projection operator corresponding to the f\/ibre at $0$
in $T^* {\mathbb R}$.
This operator, which takes $f(x)$ to the rescaled Dirac delta function $f(0) \delta(x)$, cannot be iterated.
Note that we could also consider the operator taking $f(x)$ to $r f(0) \delta(x)$ for some constant~$r$.
If the value $\delta (0)$ of the Dirac delta function at~$0$, did have a~meaning, the operator would take $\delta(x)$ to
$r \delta(0) \delta (x)$.
To have a~projection, we should take~$r$ to be $1/\delta(0)$.
If we think of~$\delta$ having an inf\/inite value at~$0$ (as in Colombeau's theory), then we could take~$r$ to be the
reciprocal of that value, an inf\/initesimal number, to get a~true projection.
This looks like some kind of renormalization; can we also interpret giving negative degrees to canonical relations as
a~similar kind of renormalization?

Note that it is possible without any technical dif\/f\/iculties to construct {\em some} projections.
For instance, say $L_1$ and $L_2$ are transversal.
Then $L_1 \times L_2$ is a~projection, even in the WW-category.
On the quantum side, the corresponding operator would, in a~simple case, map $f(x)$ to $(\int f(x) dx) \delta(x)$.
\end{itemize}

\subsection*{Acknowledgements}

Alan Weinstein would like to thank the Institut Math\'ematique de Jussieu for many years of providing a~stimulating
environment for research.
We thank Denis Auroux, Christian Blohmann, Sylvain Cappell, Alberto Cattaneo, Pavol Etingof, Theo Johnson-Freyd, Victor
Guillemin, Thomas Kragh, Jonathan Lorand, Sikimeti Mau, Pierre Schapira, Shlomo Sternberg, Katrin Wehrheim, and Chris
Woodward for helpful comments on this work.
David Li-Bland was supported by an NSF Postdoctoral Fellowship DMS-1204779; Alan Weinstein was partially supported by
NSF Grant DMS-0707137.

\pdfbookmark[1]{References}{ref}
\LastPageEnding

\end{document}